\documentclass[12pt,reqno]{amsart}
\usepackage{amsmath}
\usepackage{amsthm}
\usepackage{graphicx}
\usepackage{hyperref}
\usepackage{mathrsfs}
\usepackage{cite}
\usepackage{dsfont}
\usepackage{color}
\usepackage{bm}
\usepackage{amsfonts}
\usepackage{amssymb}
\usepackage{enumerate}
\usepackage{array}
\usepackage{cases}
\usepackage{soul}
\usepackage[margin=3cm]{geometry}
\usepackage[normalem]{ulem}

\makeatletter
\def\@settitle{\begin{center}%
		\baselineskip14\p@\relax
		\bfseries
		\uppercasenonmath\@title
		\@title
		\ifx\@subtitle\@empty\else
		\\[1ex]\uppercasenonmath\@subtitle
		\footnotesize\mdseries\@subtitle
		\fi
	\end{center}%
}
\def\subtitle#1{\gdef\@subtitle{#1}}
\def\@subtitle{}
\makeatother

\numberwithin{equation}{section}

\newtheorem{theorem}{Theorem}[section]
\newtheorem{lemma}[theorem]{Lemma}
\newtheorem{proposition}[theorem]{Proposition}

\newtheorem{definition}[theorem]{Definition}
\newtheorem{remark}[theorem]{Remark}

\definecolor{darkgreen}{rgb}{0,0.5,0}
\definecolor{darkblue}{rgb}{0,0,0.7}
\definecolor{darkred}{rgb}{0.9,0.1,0.1}
\definecolor{lightblue}{rgb}{0,0.51,1}

\newcommand{\R}{{\mathbb R}}
\newcommand{\Rd}{{\mathbb{R}^d}}
\newcommand{\N}{{\mathbb N}}
\newcommand{\E}{{\mathbb E}}
\newcommand{\Pb}{{\mathbb P}}
\newcommand{\F}{{\mathcal F}}
\newcommand{\Q}{{\mathcal Q}}
\newcommand{\T}{{\mathcal T}}

\newcommand{\dx}{D_x}
\newcommand{\dxx}{{D^2_{xx}}}
\newcommand{\dxi}{\partial_{x_i}}
\newcommand{\dxj}{\partial_{x_j}}

\newcommand{\1}{\mathds{1}}

\DeclareMathOperator{\tr}{tr}
\DeclareMathOperator{\divv}{div}

\newcommand{\xit}{X^{i,N}_t}
\newcommand{\xjt}{X^{j,N}_t}
\newcommand{\xis}{X^{i,N}_s}
\newcommand{\xjs}{X^{j,N}_s}
\newcommand{\eit}{E^{i,N}_t}
\newcommand{\ejt}{E^{j,N}_t}
\newcommand{\eis}{E^{i,N}_s}

\newcommand{\eism}{E^{i,N}_{s^-}}
\newcommand{\ejsm}{E^{j,N}_{s^-}}

\newcommand{\mun}{\mu^N_t}

\newcommand{\muns}{\mu^N_s}

\def\E{{\mathbb E}}

\def\R{{\mathbb R}}

\def\N{{\mathbb N}}

\def\P{{\mathcal P}}

\def\M{{\mathcal M}}

\def\Q{{\mathcal Q}}

\def\D{{\mathcal D}}
\def\W{{\mathcal W}}

\def\A{{\mathcal A}}
\def\F{{\mathcal F}}

\definecolor{darkspringgreen}{rgb}{0.09, 0.45, 0.27}

\title[Central Limit Theorem for a spatial stochastic epidemic]{Central Limit Theorem for a spatial stochastic epidemic model with mean field interaction}
\author{M. Hauray, E. Pardoux, Y. V. Vuong}

\begin{document}
	
	\maketitle

	\begin{abstract}
		In this article, we study an interacting particle system in the context of epidemiology where the individuals (particles) are characterized by their position and infection state. We begin with a description at the microscopic level where the displacement of individuals is driven by mean field interactions and state-dependent diffusion, whereas the epidemiological dynamic is described by the Poisson processes with an infection rate based on the distribution of other nearby individuals, also of the mean-field type. Then under suitable assumptions, a form of law of large numbers has been established to show that the associated empirical measure to the above system converges to the law of the unique solution of a nonlinear McKean-Vlasov equation. As a natural follow-up question, we study the fluctuation of this stochastic system around its limit. We prove that this fluctuation process converges to a limit process, which can be characterized as the unique solution of a linear stochastic PDE. Unlike the existing literature using a coupling approach to prove the central limit theorem for interacting particle systems, the main idea in our proof is to use a semigroup formalism and some appropriate estimates to directly study the linearized evolution equation of the fluctuation process in a suitable weighted Sobolev space and follows a Hilbertian approach.
		
	\end{abstract}


	\section{Introduction}
	
	In this paper, we study a spatial stochastic epidemic model based on the well-known SIR model, where S, I and R respectively stands for the different states of an individual. These states can vary from the compartment of Susceptible to the Infected one, and eventually to the compartment of Recovered (Removed) when the individual recovers from the illness (or dies). In fact, for many problems related to the spread of infectious disease in ecology and public health, an explicit description of spatial structure is not necessary nor advantageous. In many cases, the concept of average behavior in a large population is sufficient enough to provide the insights and extract useful information from existing data. However, the spatial component of many transmission systems is becoming increasingly important \cite{Riley-Trapman15}. Recent studies in both deterministic and stochastic epidemic models have enabled us to understand the significance of individual displacements in a population on the persistence or extinction of an endemic disease \cite{Bowong-Emakoua-Pardoux22, Emakoua22, Brauer-Driessche-Wu08, Leonard90}.
	
	\smallskip
	In our spatial model, an individual will be characterized by two features: its position and its infection state. The state $E$ varies in the set of the types $\{S,I,R\} = \{ 0,1,2\}$, where we identify $S$ with $0$, $I$ with $1$ and $R$ with $2$ in order to simplify the mathematical description. It is also useful for the representation of the jumps between the states in the epidemic dynamic. Meanwhile, the position is a continuous variable $X \in \R^d$. The addition of spatial variables complicates the standard homogeneous SIR model in two ways: by using an infection rate that depends on the distribution of surrounding population and by taking into account the individual displacements.
	
	\smallskip
	In fact, it is a natural tendency that an infected individual will infect a close neighbor more often than a more distant individual. While these different transmission behaviors are averaged in homogeneous SIR models, in our model, we use an infection rate depending on the relative distance between individuals. The infection rate between locations $x,y \in \R^d$ will be given by a function $K : \R^d \times \R^d \to  \R_+$, which is assumed to be bounded and Lipschitz. Averaging over all the infected individuals, the susceptible individual $i$ becomes infected (in other words its state jumps from $0$ to $1$) at time $t$ at the rate
	\begin{equation} \label{CLT:rate_infection}
		\frac1{N}\sum_{j=1}^N K(\xit,\xjt)\1_{\{\eit=S\}}\1_{\{\ejt=I\}}.
	\end{equation}
	
	The infectious individuals recover (in other words their state jumps from $1$ to $2$) at rate $\gamma >0$ and once an individual recovers, it is immune forever.
	
	\smallskip
	Each individual moves in $\Rd$ according to a diffusion $\sigma\bigl(\xit,\eit\bigr) dB^i_t$ which depends on both individual's state and position, and weakly interact with the other individuals in the population in the mean filed type through a kernel $V$. In this paper, the interaction kernel $V$, the diffusion strength $\sigma$ are assumed to be bounded Lipschitz continuous with respect to the position variables. Of course, this equation has a meaning on a probability space endowed with the requested Brownian motions (and Poisson point processes for the infectious-jump part of the dynamic). The Lipschitz hypothesis will be very useful to build a correct theory of existence, uniqueness to that system, and also for our results concerning the large population limit (i.e. when $N$ goes to infinity).
	
	\smallskip
	In light of the aforementioned settings, the epidemiological dynamic can be represented using Poisson point processes jumping in $\{0,1,2\}$. Now we choose a probability space $(\Omega,\F,(\F_t)_{t\geq 0},\Pb)$ equipped with $N$ independent Poisson random measures $(\Q^i)_{i=1,\dots,N}$ and $N$ Brownian motions $(B^i)_{i=1,\dots,N}$, the position and state of the individuals will evolve in time according to the following system:
	
	\begin{align}\label{CLT:interacting_system}
		\left\{
		\begin{array}{rl}
			d\xit &=  \displaystyle\frac{1}{N}\sum_{j=1}^NV \bigl(\xit,\eit,\xjt,\ejt \bigr)dt+\sigma\bigl(\xit,\eit\bigr) dB^i_t,  \\
			\eit &=   E^{i,N}_0 + \displaystyle\int_{[0,t]\times\R_+}\1_{\big\{ u \le \frac1{N} \sum_{j \neq i} K\bigl( \xis,\xjs  \bigr) \1_{(\eism,\ejsm)=(0,1)} + \gamma \1_1(\eism)\big\} } \Q^i(ds,du).
		\end{array} \right.
	\end{align}
	%
	%

	For more details concerning the origin of this model and its interest, we refer readers to the previous paper of the authors \cite{Hauray-Pardoux-Vuong22}.
	
	\smallskip
	
	In the study of a system composed of $N$ particles, one of the most important objects is the empirical measure that can help us fully describe the whole dynamic. In this paper, let us introduce the empirical measure process associated to the above system consisting of $N$ individuals $\bigl(\xit,\eit)$, $i=1,\ldots,N$ defined by
	\[
	t  \mapsto \mun=\frac{1}{N}\sum_{i=1}^N\delta_{(\xit,\eit)},
	\]
	where $\delta_{(x,e)}$ is the Dirac measure at point $(x,e)\in \R^d\times \{0,1,2\}$.
	
	\smallskip
	Under suitable assumptions in  \cite{Hauray-Pardoux-Vuong22}, we established a conditional propagation of chaos result (in the presence of a common noise $\sigma^0(X^{i,N}_t,E^{i,N}_t)dB^0_t$), which states that conditionally to the common noise, the individuals are asymptotically independent and the stochastic dynamic converges to a random nonlinear McKean-Vlasov process when the population size tends to infinity. And as a consequence, the associated empirical measure converges to the unique solution of a stochastic mean-field PDE driven by the common noise.
	
	\smallskip
	In this work, we only treat the case without the common noise. As a special case of the results obtained in \cite{Hauray-Pardoux-Vuong22} (with $\sigma^0=0$), we can also show that when $N\to \infty$, the empirical measure $\mu^N_t$ converges to $\mu_t$ the law of the unique solution to the following nonlinear McKean-Vlasov equation
	\begin{equation}\label{CLT:limit_SDE}
		\begin{cases}
			dX_t &= V_{\mu_t}(X_t,E_t) dt+\sigma(X_t,E_t)dB_t,\\
			E_t &= E_0 + \displaystyle\int_{[0,t]\times\R_+}\1_{\big\{u\leq K_{\mu_s}(X_s)\1_0(E_{s^-})+\gamma \1_1(E_{s^-})\big\}}\Q(ds,du),\\
			\mu_t &= \mathcal{L}\left(X_t, E_t\right).
		\end{cases}
	\end{equation}
	
	As typical with McKean-Vlasov dynamics, the limit measure $\mu_t$ can also be characterized as the unique solution of a nonlinear partial differential equation. That PDE is called the forward Kolmogorov equation associated to the McKean-Vlasov SDE \eqref{CLT:limit_SDE} and given by the following equation
	\begin{equation}\label{limit_SPDE}
		\begin{aligned}
			d\mu_t    =&  -\dx\cdot  \left(V_{\mu_t}\mu_t\right)dt +\frac{1}{2}\tr\left[\dxx \big((\sigma\sigma^T) \mu_t\big)\right]dt\\
			&+K_{\mu_t}\bigl(\1_{e=1} - \1_{e=0} \bigr)\mu_t(dx,0) dt+\gamma \bigl(\1_{e=2} - \1_{e=1} \bigr)\mu_t(dx,1) dt.
		\end{aligned}
	\end{equation}
	
	Now, as a natural follow-up question after studying the law of large numbers, the aim of this paper is to look for a limit theorem for the fluctuation process of $\mu^N_t$ around its limit $\mu_t$.
	
	\smallskip
	In the previous work, a quantitative law of large numbers is established in the Wasserstein distance, which roughly shows that
	\begin{equation}\label{CLT:rate_of_convergence}
		\E\left[W_1\left(\mu^N_t ,\mu_t \right)\right]\le C(t)
		\begin{cases}
			N^{-1/2},& d =1, 
			\\ N^{-1/2}\log N,& d=2, 
			\\ N^{-1/d},& d \ge 3. 
		\end{cases}
	\end{equation}
	Moreover, it is well-known that the $1$-Wasserstein distance used in~\eqref{CLT:rate_of_convergence} is equivalent to its dual formulation,
	\begin{equation}\label{eq:Kontorovich-Rubistein-dual-formula}
		W_1(\mu^N_t,\mu_t)=\sup\left\{\int_{\R^d\times\{0,1,2\}}\phi\left(\mu^N_t-\mu_t\right)\mid\; \phi : \R^d\times\{0,1,2\} \to \R \text{ with } Lip(\phi)\leq 1\right\},
	\end{equation}
	which apparently shows the strong dependence of the rate of convergence on the regularity of test functions. Therefore, to recover the right order of $N^{1/2}$-normalization  as the classical central limit theorem, we need to modify the regularity of test functions. This point will be further clarified in the next section.
	
	\smallskip
	Now we consider the following fluctuation process with the $N^{1/2}$ scaling: $$\eta^N_t=\sqrt{N}\bigl(\mu^N_t-\mu_t\bigr), \quad t\in [0,T].$$
	
	\smallskip
	Following the Hilbertian approach used in \cite{Fernandez-Meleard97,Meleard96,Jourdain-Meleard98}, we can prove a central limit theorem for the sequence of the fluctuation processes $(\eta^N)_{N\ge 1}$ in an appropriate space of distributions. The limit process of the normalized fluctuation processes can be described as the unique solution a linear stochastic partial differential equation driven by space-time white noises. In order to achieve this, we regard the fluctuation process $\eta^N_t$ as a process taking values in a Hilbert space, which we consider as the dual of some Sobolev space of test functions. The regularity of that dual space corresponding to the regularity of test functions will be decided by the martingale term appearing in the evolution equation of the fluctuations as well as the form of generators in that equation.
	
	\smallskip
	It is worth noticing that the Sobolev spaces used in the present paper are not exactly the classical one and they must be refined. Indeed, we will study a class of the weighted Sobolev spaces with polynomial weights, see the definition in subsection~\ref{sec:weighted_sobolev_spaces}. The importance of the weight will be explained in the proof, provided that the weight satisfies some suitable integrability properties. Moreover, we observe that the dimension $d$ plays a crucial role in the rate of convergence~\eqref{CLT:rate_of_convergence} and it is also well-known that the Sobolev embeddings depend strongly on the dimension of the space. This will help us identify the right level of smoothness.
	
	\smallskip
	Let us now discuss the main differences between our results and the previous one in the exiting literature. In fact, this kind of spatial epidemic model have been studied by Emakoua et al.~\cite{Bowong-Emakoua-Pardoux22,Emakoua22} with the same SIR epidemic dynamic but with a simpler model for the displacement of individuals (individual's movements follow independent Brownian motions on a compact torus in \cite{Bowong-Emakoua-Pardoux22}, and follow independent diffusion processes in \cite{Emakoua22}), where the mean field interactions between individuals through the kernel $V$ are not taken into account. This leads to the main difficulty in comparison with the previous works due to the presence of nonlocal terms in the evolution equation of the fluctuations $(\eta^N)_{N\ge 1}$. In contrast to the independence of individual's movements in \cite{Bowong-Emakoua-Pardoux22, Emakoua22}, these nonlocal terms are created by mean field interactions and they do not allow to obtain directly good estimates for the norm of fluctuations in the weighted Sobolev spaces.
	
	\smallskip
	The Hilbertian approach used in this work has already been used to prove central limit theorem in the context of interacting particle systems \cite{Meleard96, Fernandez-Meleard97, Jourdain-Meleard98}, mean field games \cite{Delarue-Lacker-Ramanan19}, mean field age-dependent Hawkes processes \cite{Chevallier17}, neuron networks \cite{Sirignano-Spiliopoulos20}. In \cite{Meleard96, Fernandez-Meleard97}, Méléard et al. developed a coupling method used in \cite{Sznitman91} and \cite{Hitsuda-Mitoma86} with some relaxations on the initial conditions and coefficients. The authors provide a sharper estimate on the control of the couplings (instead of the original one of order $2$ in the proof of the quantitative law of large numbers),
	\begin{equation}\label{sharper_estimate_order_4}
		\E\left[\sup_{t\le T}\big|X^{i,N}_t-X^i_t\big|^4\right]\le \frac{C_T}{N^2},
	\end{equation}
	where $(X^i),\; i=\{1,\ldots, N\}$ are i.i.d copies of the unique solution to the limit SDE of their original system. This estimate of order 4 requires a careful computation of the covariance between the pairs $\big(X^{i,N},X^i\big)$ and takes advantage of the independence between the particles $(X^i),\; i=\{1,\ldots, N\}$.
	
	\smallskip
	In \cite{Delarue-Lacker-Ramanan19}, Delarue et al. also used the coupling method and this estimate of order 4 to prove central limit theorem for a system consisting of $N$ agents in the context of mean field games. The main idea is to use the solution to the mean field limit to construct an associated McKean-Vlasov interacting system of $N$ particles that is suﬃciently close to the original system for large $N$, then derive the central limit theorem for the latter from the central limit theorem for the former.
	
	\smallskip
	However, the main reason that prevents us from applying this coupling method to prove central limit theorem is that the authors in the aforementioned articles only work in a continuous framework and rely strongly on the estimate of order 4 \eqref{sharper_estimate_order_4}.  In contrast, the individuals in our model possess both continuous and discrete features. In the previous work, we have pointed out the compulsion of using estimates of order~$1$ for the control of the couplings (see Remark $3.2$ in \cite{Hauray-Pardoux-Vuong22}). As usual when working with jump processes, we can not get higher rate for the moment estimates as in \eqref{sharper_estimate_order_4}.  Hence the standard trick used for diﬀusion processes is useless in this case. To solve this difficulty, the author in \cite{Chevallier17} developed the above coupling method for a specific mean field interacting age-dependent Hawkes process. A refined version of the higher order estimates \eqref{sharper_estimate_order_4} is provided by estimating the coupling in the total variation sense.
	
	\smallskip
	Unlike the articles listed above, where the coupling method is used to prove the central limit theorem, in the proof of the present paper, we use the semigroup formalism and some appropriate estimates to directly study the linearized evolution equation of the fluctuation process in a suitable weighted Sobolev space.
	It will be shown that under some suitable assumptions on the initial conditions and the smoothness of the coefficients, the fluctuation processes $(\eta^N_t)_{N\ge 1}$ belong uniformly in $N$ and $t$ to the weighted Sobolev spaces $H^{-(1+D), 2 D}$ and $H^{-(2+2 D), D}$ (see the definition in subsection~\ref{sec:weighted_sobolev_spaces}, with $D:=\lceil d/2\rceil$ ). Then we prove the tightness of the pre-limit fluctuation process in $\D\big([0,T], H^{-(2+2 D), D}\big)$ by using appropriate compact embeddings. We also show that the Hilbert space $H^{-(2+2 D), D}$ is sharp in some sense: it has the smallest regularity order as possible in the class of Sobolev spaces with polynomial weights where we can obtain the tightness result. Finally, we complete the proof of convergence of the sequence $(\eta^N)_{N\ge 1}$ by identifying the limit fluctuation process $\eta$ as the unique solution of a linear stochastic partial differential equation.
	

	\smallskip
	\paragraph{\bf Organisation of the paper.}
	In Section~\ref{sec:preliminaries_and_mainresult}, we provide some preliminaries on the weighted Sobolev spaces and state the main results. Section~\ref{sec:tightness} is devoted to prove the tightness of the pre-limit fluctuation process and the martingale terms appearing in the evolution equation. In order to do this, we first establish some key estimates in dual Sobolev norms and then take advantage of the Hilbert structure of the Sobolev spaces to prove the tightness results. Section~\ref{sec:characterzation_of_the_limit} contains the proof of the main Theorem~\ref{thm:main_theorem_etaN}, and we give a characterization of the limit fluctuation process as the unique solution to a linear SPDE driven by space-time white noises.

	\section{Preliminaries and main result}\label{sec:preliminaries_and_mainresult}
	\subsection{Preliminaries on weighted Sobolev spaces}\label{sec:weighted_sobolev_spaces}
	This section is devoted to the definitions and some technical results related to the Sobolev spaces with polynomial weights used in this paper. This kind of weighted Sobolev spaces was first introduced in \cite{Metivier84}, see also \cite{Fernandez-Meleard97}.
	
	\medskip
	\paragraph{\textit{Weighted Sobolev spaces.}} For all $j\in\N,\, \alpha >0,\, g\in C^j(\R^d)$, we define
	$$\|g\|_{j,\alpha}:=\left(\sum_{|k|\leq j} \int_{\R^d} \frac{|D^kg(x)|^2}{(1+|x|^2)^{\alpha}}dx\right)^{1/2},
	$$
	where $k=(k_1,\ldots, k_d), \; |k|=k_1+\cdots+k_d$.
	
	\smallskip
	Let $H^{j,\alpha}(\R^d)$ be the completion of the space consisting of all functions $g\in C^{\infty}(\R^d)$ with compact support with respect to the $\|\cdot\|_{j,\alpha}$ norm. $H^{j,\alpha}$ equipped with this norm is a Hilbert space. We denote by $H^{-j,\alpha}$ its dual space.
	
	\smallskip
	Let $C^{j, \alpha}$ be the space of functions $g$ with continuous partial derivatives up to order $j$ and satisfies 
	\[
	\lim_{|x| \rightarrow \infty} \frac{\left|D^{k} g(x)\right|}{1+|x|^{\alpha}}=0, \quad \forall\, |k| \le j.
	\]
	
	\smallskip
	This space is normed with
	$$
	\|g\|_{C^{j,\alpha}}=\sum_{|k| \le j} \sup _{x \in \mathbb{R}^{d}} \frac{\left|D^{k} g(x)\right|}{1+|x|^{\alpha}}.
	$$
	
	Noticing that for all $j\ge 0$, $C^{j, 0}\equiv C^j_b$, the space of $C^j$ functions with bounded derivatives of all order less than $j$.
	
	\medskip
	\paragraph{\textit{Sobolev embeddings.}} We recall the some continuous embeddings related to the Sobolev spaces defined above, which are useful in some proofs in the rest of this paper. For more details, see e.g. \cite{Adams03}, \cite{Fernandez-Meleard97}.
	
	\smallskip
	We have
	\begin{align}
		C^j_b&\hookrightarrow H^{j,\alpha}, & &j\ge 0, \quad \alpha >d/2,\;\Big(\text{so that } \int_{\R^d} 1/(1+|x|^{2\alpha}) dx<+\infty\Big), \label{embedding:C_to_W}\\
		H^{j+m, \alpha} &\hookrightarrow C^{j, \alpha}, & &j\ge0, \quad m>d/2,\quad \alpha \ge 0, \label{embedding:W_to_C}
	\end{align}
	i.e. there exists $C_1,\, C_2$ (that depends on $m, j$ and $\alpha$) such that
	\begin{align*}
		\big\|g\big\|_{H^{j, \alpha}} &\le C_1\big\|g\big\|_{C^j_b}, \\
		\big\|g\big\|_{C^{j, \alpha}} &\le C_2\big\|g\big\|_{H^{j+m, \alpha}}.
	\end{align*}
	
	\smallskip
	Moreover, using the embedding \eqref{embedding:W_to_C}, we can prove that
	\begin{align}\label{embedding:W_to_W}
		H^{j+m, \alpha} &\hookrightarrow_{c} H^{j, \alpha+\beta}, & &j \ge 0, \quad m>d/2, \quad \alpha \ge 0, \quad \beta>d/2,
	\end{align}
	where $\hookrightarrow_{c}$ means that the embedding is compact.
	
	\smallskip
	We also deduce the following dual embeddings:
	\begin{align}
		H^{-j,\alpha} &\hookrightarrow C^{-j,0}, & &j\ge 0, \quad \alpha > d/2, \label{embedding:Wdual_to_Cdual}\\
		C^{-j,\alpha} &\hookrightarrow H^{-(j+m), \alpha}, & & j\ge0, \quad m>d/2,\quad \alpha \ge 0, \label{embedding:Cdual_to_Wdual}\\
		H^{-j, \alpha+\beta} &\hookrightarrow_{c} H^{-(j+m), \alpha}, & & j\ge0, \quad m>d/2, \quad  \alpha \ge 0, \quad \beta>d / 2 \label{embedding:Wdual_to_Wdual}.
	\end{align}

	\medskip
	\paragraph{\textit{Hilbert structures.}} In the next sections, once $(\phi_{k})_{k\geq 1}$ is mentioned, it always denotes an orthonormal basis of $H^{j,\alpha}$  composed of $C^{\infty}$ functions with compact support. The existence of this basis follows from the fact that the functions of class $C^{\infty}$ with compact support are dense in $H^{j,\alpha}$. Moreover, if $(\phi_{k})_{k\geq 1}$ is an orthonormal basis of $H^{j,\alpha}(\R^d)$ and $u$ belongs to $H^{-j,\alpha}(\R^d)$ then Parseval's identity give us the following representation
	\begin{equation}\label{Parseval_identity}
		\big\|u\big\|^{2}_{-j,\alpha}=\sum_{k\geq 1} \left< u,\phi_{k}\right>^{2}
	\end{equation}
	
	We also note that for any distribution $\mu=(\mu^0,\mu^1,\mu^2) \in \D(\R^d\times\{0,1,2\})$, we can define
	\begin{equation*}
		\big\|\mu\big\|^{2}_{H^{-j,\alpha}(\R^d\times\{0,1,2\})}:= \big\|\mu^0\big\|^{2}_{H^{-j,\alpha}(\R^d)}+\big\|\mu^1\big\|^{2}_{H^{-j,\alpha}(\R^d)}+\big\|\mu^2\big\|^{2}_{H^{-j,\alpha}(\R^d)}.
	\end{equation*}
	For any test function $\phi=(\phi^0,\phi^1,\phi^2)$ and $\mu=(\mu^0,\mu^1,\mu^2)$,
	\begin{equation*}
		\big<\mu,\phi\big>=\sum_{i=0}^{2} \big<\mu^i,\phi^i\big>=\sum_{i=0}^{2} \int_{\R^d} \phi^i(x)\mu^i(dx).
	\end{equation*}

	\subsection{Main results}
	In this section, we rigorously describe the evolution equation of the fluctuation process and state the main results. As in the previous paper \cite{Hauray-Pardoux-Vuong22}, by using Itô's formula we showed that the evolution of the empirical measure process $\mu^N_t$ satisfies the following equation:
	\begin{equation}\label{evolution_empirical_measure}
		\begin{aligned}
			\left \langle \mun,\phi \right \rangle=&\left\langle\mu^N_0,\phi \right\rangle+\int_{0}^{t}\left\langle\muns, \dx\phi\cdot V_{\muns}\right\rangle ds +\frac{1}{2}\int_{0}^{t}\left\langle\muns,\tr\left[(\sigma\sigma^T)\dxx\phi\right]\right\rangle ds\\
			&+\int_0^t\left\langle \mu^N_s(dx,0),K_{\mu^N_s}(\1_{e=1}-\1_{e=0})\phi\right\rangle ds+\int_0^t\left\langle \mu^N_s(dx,1),\gamma (\1_{e=2}-\1_{e=1})\phi \right\rangle ds\\
			&+M^N_t(\phi) ,
		\end{aligned}
	\end{equation}
	where $M^N_t(\phi)$ is a martingale which converges to $0$,
	\begin{align*}
		M^N_t(\phi) =&\frac{1}{N}\sum_{i=1}^N\int_{0}^{t}\dx\phi(\xis,\eis)\sigma(\xis,\eis)dB^i_s\\
		&+\frac{1}{N}\sum_{i=1}^N\int_{[0,t]\times\R_+}\big(\phi(X^{i,N}_s,E^{i,N}_s)-\phi(X^{i,N}_s,E^{i,N}_{s^-})\big)\times\\
		&\hspace{3cm}\times\1_{\big\{u\leq K_{\mu^N_s}(X^{i,N}_s) \1_0(E^{i,N}_{s^-}) + \gamma \1_1(E^{i,N}_{s^-}) \big\}}\bar{Q}^i(ds,du).
	\end{align*}

	Subtracting the equation \eqref{limit_SPDE} to the evolution equation \eqref{evolution_empirical_measure} of the empirical measure $\mu^N_t$, and then multiplying by $N^{1/2}$, we can obtain the evolution equation of the fluctuation process $\eta^N_t$ as the following:
	\begin{equation}\label{evolution_fluctuation_process}
		\begin{aligned}
			\left \langle \eta^N_t,\phi \right \rangle=&\left\langle\eta^N_0,\phi \right\rangle+\frac{1}{2}\int_{0}^{t}\left\langle\eta^N_s,\tr\left[(\sigma\sigma^T)\dxx\phi\right]\right\rangle ds\\
			&+\int_{0}^{t}\left\langle\eta^N_s, \dx\phi\cdot V_{\mu^N_s}\right\rangle ds+\int_{0}^{t}\left\langle\mu_s, \dx\phi\cdot V_{\eta^N_s}\right\rangle ds\\
			&+\int_0^t\left\langle\eta^N_s(dx,0), K_{\mu^N_s}(\1_{e=1}-\1_{e=0})\phi\right\rangle ds+\int_0^t\left\langle \mu_s(dx,0),K_{\eta^N_s}(\1_{e=1}-\1_{e=0})\phi\right\rangle ds\\
			&+\int_0^t\left\langle \eta^N_s(dx,1),\gamma (\1_{e=2}-\1_{e=1})\phi \right\rangle ds +\sqrt{N}M^N_t(\phi).
		\end{aligned}
	\end{equation}
	
	It is worth noting that the two terms in the second line on the r.h.s. are created by linearizing the nonlinear term $\left\langle\muns, \dx\phi\cdot V_{\muns}\right\rangle$, whereas the two terms in the third line are the linearization of $\left\langle \mu^N_s(dx,0),K_{\mu^N_s}(\1_{e=1}-\1_{e=0})\phi\right\rangle$. In contrast to the law of large numbers, the martingale term in \eqref{evolution_fluctuation_process} does not go to $0$ when $N$ tends to infinity. Instead of vanishing, the renormalized martingale $\sqrt{N}M^N_t$ is expected to converge to some Gaussian process.
	
	\smallskip
	Before giving a statement about the convergence, the first problem one needs to overcome is to find a suitable space in which both $\eta^N$ and its limit belong. We want to prove that the fluctuation $\eta_{t}^N$ belongs to some weighted Sobolev space $H^{-j,\alpha}$ uniformly in $N$ and $t\in [0,T]$. By taking an orthonormal basis $(\phi_{k})_{k\geq 1}$ of the Sobolev space $H^{j,\alpha}$ as in \eqref{Parseval_identity}, our desire is to get the following
	\begin{equation}\label{eq:uniform_estimate_desired}
		\sup_{N\ge 1} \sup_{t\le T}\E\left[\sum_{k\ge 1}\left\langle\eta_{t}^N, \phi_{k}\right\rangle^{2}\right]=\sup_{N\ge 1} \sup _{t \le T}\E\left[\big\|\eta_{t}^N\big\|_{-j,\alpha}^{2}\right]<+\infty.
	\end{equation}
	
	To see the impact of the regularity of test functions on the estimates of the fluctuation process $\eta^N_t$ in the dual spaces, let us provide in the following a simple example on the class of functions with bounded Lipschitz constant, where we can compute properly by using the Kantorovich-Rubistein duality \eqref{eq:Kontorovich-Rubistein-dual-formula}. Indeed, from the quantitative law of large number, we have
	\begin{align*}
		\E\left[\sup_{\phi\in Lip(1)}\left|\left\langle\eta_{t}^N, \phi\right\rangle\right|\right]=&\E\left[\sup_{\phi\in Lip(1)}\left|\left\langle\sqrt{N}(\mu^N_t-\mu_t), \phi\right\rangle\right|\right]\\
		=&\sqrt{N}\E\left[W_1(\mu^N_t,\mu_t)\right]\\
		\le& C(t)\left(\sqrt{N}\E\left[W_1\left(\mu^N_0,\mu_0\right)\right]+\begin{cases}
			1,& d =1, 
			\\ \log N,& d=2, 
			\\ N^{(d-2)/2d},& d \ge 3. 
		\end{cases}\right).
	\end{align*}
	Since $\big\{(X^{i,N}_0,E^{i,N}_0)\big\}_{N\ge 1}$ are i.i.d. with the initial law $\mu_0$, the classical central limit theorem ensures at initial time that $\eta^N_0$ converges in law to a limit $\eta_0$, which is a Gaussian. However, the above estimate is obviously not enough to guarantee central limit theorem for the fluctuation process when it evolves in time, and even the uniform estimate \eqref{eq:uniform_estimate_desired} fails when the dimension $d$ is large. Therefore, in order to obtain the needed estimates and recover the right order for convergence in central limit theorem, test functions indeed must be more regular than only Lipschitz.
	
	\medskip
	Before stating the main results, let us introduce the assumptions made for the initial condition and the coefficients throughout this paper. 
	
	\paragraph{Assumptions.}

	We fix $D:=\lceil d/2\rceil$.
	
	\smallskip
	Assumption \textbf{$\A_1$}. \label{assumption_A1}
	$\sup_{N\ge 1}\max_{1\le i\le N}\E\left[\big|X_0^{i,N}\big|^{4D}\right]<+\infty$.
	
	\smallskip
	Assumption \textbf{$\A_2$}. \label{assumption_A2}
	The functions $V, K$ belong to class $C^{1+D}_b$ and $\sigma \in C^{4D+5}_b$. We also assume that the symmetric matrix $[\sigma\sigma^{\dag}]$ is uniformly positive definite.
	
	\smallskip
	Assumption \textbf{$\A_3$}. \label{assumption_A3}
	The functions $V, K$ belong to class $C^{2+2D}_b$.
	
	\medskip
	The assumptions~\textbf{$\A_1,\, \A_2$} above are essential to prove the propagation of moments and the tightness results in the next section. We also notice that with the hypothesis $[\sigma\sigma^{\dag}]$ is uniformly positive definite, the operator $A:=\frac{1}{2}\dx\cdot\left(\sigma{\sigma}^{\dag}\dx\right)$ is uniformly elliptic, i.e.
	\[\sum_{i,j=1}^{d} (\sigma\sigma^{\dag})_{ij}(x,e)\xi_i\xi_j\ge \lambda |\xi|^2, \quad \forall (x,e)\in\R^d\times\{0,1,2\}, \;\xi \in \R^d, \]
	for some positive constant $\lambda$. This assumption allows us to perform some crucial estimates in the proof of Proposition~\ref{prop:uniform_estimate_etaN}. In order to characterize the limit fluctuation process as in the statement of the central limit theorem~\ref{thm:main_theorem_etaN}, more regularity on the coefficients will be required and given in Assumption~\textbf{$\A_3$}.
	
	\medskip
	It is shown in the following that under appropriate assumptions on the initial conditions and the smoothness of the coefficients, the fluctuation processes $(\eta^N_t)_{N\ge 1}$ belong uniformly in $N$ and $t$ to $H^{-(1+D), 2 D}$. 
	
	\begin{proposition}\label{prop:uniform_estimate_etaN}
		Let $T>0$. Under Assumptions \textbf{$\A_1,\, \A_2$}, the fluctuation process $(\eta^{N}_t)_{t\le T}$ belongs to $H^{-(1+D), 2D}$ uniformly in $t$ and $N$, i.e.
		\begin{equation}\label{eq:uniform_estimate_etaN}
			\sup_{N\ge 1} \E\left[\sup _{t \le T}\big\|\eta^{N}_t\big\|_{-(1+D),2D}^{2}\right]<+\infty.
		\end{equation}
	\end{proposition}
	
	Then we prove the tightness of the fluctuation process $\eta^N_t$ in $\D\big([0,T], H^{-(2+2 D), D}\big)$ by using the embeddings described in subsection~\ref{sec:weighted_sobolev_spaces}.
	
	\begin{proposition}\label{prop:tightness_of_eta^N}
		Under Assumptions \textbf{$\A_1,\, \A_2$}, the sequence of the laws of $(\eta^{N})_{N \ge 1}$ is tight in $\D\left([0, T], H^{-(2+2 D), D}\right)$.
	\end{proposition}	
	
	The main result of this paper will be stated below. It identifies the limit fluctuation process $\eta$ as the unique solution of a linear stochastic partial differential equation.
	
	\begin{theorem}\label{thm:main_theorem_etaN}
		Under Assumptions \textbf{$\A_1, \A_2, \A_3$}, the sequence of fluctuation processes $\bigl(\eta^{N}\bigr)_{N \geq 1}$ converges in law in $\D\left([0, T], H^{-(2+2 D), D}\right)$ to a process $\eta$ which solves the following equation

		\begin{equation}\label{eq:identifying_limit_eta}
			\begin{aligned}
				\eta_t=&\eta_0+\frac{1}{2}\int_{0}^{t}\tr\left[\dxx\big((\sigma\sigma^T)\eta_s\big)\right]ds -\int_{0}^{t}\dx\cdot(\eta_sV_{\mu_s})ds-\int_{0}^{t}\dx\cdot(\mu_s V_{\eta_s})ds\\
				&+\int_0^t K_{\mu_s}(\1_{e=1}-\1_{e=0})\eta_s(dx,0)ds+\int_0^t K_{\eta_s}(\1_{e=1}-\1_{e=0})\mu_s(dx,0)ds\\
				&+\int_0^t \gamma(\1_{e=2}-\1_{e=1})\eta_s(dx,1)ds +\W_t,
			\end{aligned}
		\end{equation}
		where $\W_t$ is a continuous centered Gaussian process with values in $H^{-(2+2 D), D}$ and covariance is given by: For all $\phi_1, \phi_2 \in H^{(2+2 D), D}$, for any $s, t \in [0,T]$,
		\begin{equation}\label{eq:covariance_Gaussian_identifying_limit_eta}
			\begin{aligned}
				\E\left[\W_t(\phi_1)\W_s(\phi_2)\right]=&\int_0^{t\wedge s} \left<\mu_r, \sigma\sigma^T \dx\phi_1\cdot\dx\phi_2  \right>dr\\
				&+\int_0^{t\wedge s} \big\langle \mu_s(dx,0),  K_{\mu_r(dx,1)}\phi_1\phi_2\big\rangle dr + \int_0^{t\wedge s} \big\langle \mu_r(dx,1), \gamma\phi_1\phi_2\big\rangle dr.
			\end{aligned}
		\end{equation}
		
	\end{theorem}

	\section{Tightness}\label{sec:tightness}
	
	\subsection{Preliminary estimates}\label{sec:CLT:preliminary_key_estimates}
	In this section, we first prove some useful estimates which are the technical steps in the proof of tightness and convergence in the next sections.
	
	\medskip
	
	We first recall a fundamental result which states that the initial condition \textbf{$\A_1$} propagates finite moments uniformly in $N$ and time $t\in[0,T]$. The proof of this result is classical.
	\begin{lemma}\label{lem:finite_moment_estimates}
		For any $T>0$, there exists a constant $C_T$ such that
		\begin{align*}
			\sup_{N\ge 1}\E\bigg[\sup _{t \le T}\big|X_{t}^{i, N}\big|^{4D}\bigg]&\le C_T, \quad \forall \,1 \le i \le N,\\
			\E\bigg[\sup _{t \le T}\left|X_{t}\right|^{4D}\bigg]&\le C_T.
		\end{align*}
	\end{lemma}
	
	\begin{remark}
		By the definition of the empirical measure $\mu^{N}_t$ and its limit $\mu_t$, we can easily deduce from Lemma \ref{lem:finite_moment_estimates} that
		\begin{align*}
			\sup_{N\ge 1}\E\bigg[\sup _{t \le T}\left<\mu^{N}_t, |\cdot|^{4D}\right>\bigg]&\le C_T, \\
			\E\bigg[\sup _{t \le T}\left<\mu_t, |\cdot|^{4D}\right>\bigg]&\le C_T.
		\end{align*}
	\end{remark}

	\medskip
	Next, we give some useful estimates of several linear operators on $H^{j,\alpha}$. We may use them many times in the next sections.
	
	\begin{lemma}\label{lem:control_linear_forms}
		For any fixed $\alpha\ge 0$, $j\ge 1+D$ and $x,y\in \R^d$, the mappings $\delta_{x},\Lambda_{x,y},\Psi _x:H^{j,\alpha}\to \mathbb{R}$, defined by \[\delta_{x}(\phi):=\phi(x);\quad \Lambda_{x,y}(\phi):=\phi(x)-\phi(y);\quad \Psi _x(\phi):=(\divv\phi)(x)\] are continuous linear forms, and we have
		\begin{equation}
			\begin{aligned}
				\big\| \delta_{x} \big\|_{-j,\alpha}\leq& K (1+ |x|^{\alpha}),\\
				\big\| \Lambda_{x,y} \big\|_{-j,\alpha}\leq& K (1+ |x|^{\alpha}+|y|^{\alpha}),\\
				\big\| \Psi _x\big\|_{-j,\alpha}\leq& K (1+ |x|^{\alpha}).
			\end{aligned}
		\end{equation}
	\end{lemma}
	\begin{proof}
		We prove the first estimate by applying the embedding \eqref{embedding:W_to_C},
		\begin{equation}\label{eq:proof:linear_forms_ineq1}
			|\delta_{x}(\phi)|=|\phi(x)|\leq \big\|\phi\big\|_{C^{0,\alpha}} (1+|x|^{\alpha})\leq K \big\|\phi\big\|_{j,\alpha} (1+|x|^{\alpha}),\quad \forall\, j\ge D, \alpha \ge 0.
		\end{equation}
		
		Using the definition of dual norms of linear mappings, we have
		\[\big\| \delta_{x} \big\|_{-j,\alpha}=\sup_{0\ne\phi\in H^{j,\alpha}}\frac{|\delta_{x}(\phi)|}{\big\|\phi\big\|_{j,\alpha}} \leq K (1+ |x|^{\alpha}).\]
		
		The estimate for $\Lambda_{x,y}$ follows \eqref{eq:proof:linear_forms_ineq1} since \[|\Lambda_{x,y}(\phi)|\leq |\phi(x)|+|\phi(y)|= |\delta_{x}(\phi)| + |\delta_{y}(\phi)|.\]
		
		A similar argument holds true for $\Psi _x$ with $j\ge D+1$ and $\alpha \ge 0$.
		
	\end{proof}
	
	\subsection{Decomposition of the fluctuations}\label{section:Decomposition_fluc}
	
	In this section, we will describe the fluctuation process $(\eta^N_t)_{t\ge 0}$ explicitly in terms of each epidemiological state S, I and R. On the one hand, this turns the equation \eqref{evolution_fluctuation_process} to a system consisting of three equations. On the other hand, rewriting the evolution equation of fluctuation process as a system seems to be compatible with our strategy to prove the convergence in the next section. Indeed, we will represent the linearized equation~\eqref{evolution_fluctuation_process} in a semigroup formalism, and take advantage of some useful estimates in the semigroup approach to prove the key estimate~\eqref{eq:uniform_estimate_etaN}. For that reason, in order to make the semigroup representation of the evolution equation \eqref{evolution_fluctuation_process} less complex, we consider its projections on $\M(\R^d)$ for each epidemiological state separately. For more details concerning this semigroup representation, see Section \ref{sec:proof_main_Prop-(1)}.
	
	\medskip
	For $x\in \R^d$, let
	\begin{align*}
		\mu^{S,N}(x):=\mu^N(x,0),\\
		\mu^{I,N}(x):=\mu^N(x,1),\\
		\mu^{R,N}(x):=\mu^N(x,2).
	\end{align*}
	We regard $\mu^{S,N},\mu^{I,N},\mu^{R,N}$ as \text{càdlàg} processes taking values in the space of finite measures on $\R^d$, equipped with the Skorohod topology.
	
	\smallskip
	For each $e\in{\{S,I,R\}}$, we introduce the following alternative notations
	\begin{align*}
		\sigma^e(\cdot):=&\sigma(\cdot,e),\\
		V^e_\mu(\cdot):=&V_\mu(\cdot,e)=\big\langle V(\cdot,e,y,f), \mu(dy,df)\big\rangle
	\end{align*}
	to adapt with the measures on $\R^d$.
	
	\smallskip
	We also note that somewhere the notion $V^e$ will be assigned to a function of there variables on $\R^d\times\R^d\times\{0,1,2\}$, $V^e(\cdot,\cdot,\cdot):=V(\cdot,e,\cdot,\cdot)$.
	
	\smallskip
	Now as usual, by using Itô's formula we can derive the evolution equation for the empirical measures $\mu^{S,N},\mu^{I,N},\mu^{R,N}$. Indeed, for any test function $\phi\in C^2_b(\R^d)$, we have the following system which is equivalent to equation \eqref{evolution_empirical_measure}:
	\begin{align}
		\big\langle \mu^{S,N}_t,\phi \big\rangle
		=&\big\langle\mu^{S,N}_0,\phi \big\rangle+\frac{1}{2}\int_{0}^{t}\big\langle\mu^{S,N}_s,\tr\big[(\sigma^S{\sigma^S}^{\dag})\dxx\phi\big]\big\rangle ds + \int_{0}^{t}\big\langle\mu^{S,N}_s, \dx\phi\cdot V^S_{\mu^N_s}\big\rangle ds \nonumber\\
		&-\int_0^t\big\langle\mu^{S,N}_s, \phi K_{\mu^{I,N}_s}\big\rangle ds+ M^{S,N}_t(\phi), \label{eq:empirical_mu^SN}\\
		\big\langle \mu^{I,N}_t,\phi \big\rangle
		=&\big\langle\mu^{I,N}_0,\phi \big\rangle+\frac{1}{2}\int_{0}^{t}\big\langle\mu^{I,N}_s,\tr\big[(\sigma^I{\sigma^I}^{\dag})\dxx\phi\big]\big\rangle ds + \int_{0}^{t}\big\langle\mu^{I,N}_s, \dx\phi\cdot V^I_{\mu^N_s}\big\rangle ds \nonumber\\
		&+\int_0^t\big\langle\mu^{S,N}_s, \phi K_{\mu^{I,N}_s}\big\rangle ds-\gamma\int_0^t\big\langle\mu^{I,N}_s,\phi \big\rangle ds+ M^{I,N}_t(\phi), \label{eq:empirical_mu^IN}\\
		\big\langle \mu^{R,N}_t,\phi \big\rangle
		=&\big\langle\mu^{R,N}_0,\phi \big\rangle+\frac{1}{2}\int_{0}^{t}\big\langle\mu^{R,N}_s,\tr\big[(\sigma^R{\sigma^R}^{\dag})\dxx\phi\big]\big\rangle ds + \int_{0}^{t}\big\langle\mu^{R,N}_s, \dx\phi\cdot V^R_{\mu^N_s}\big\rangle ds \nonumber\\
		&+\gamma\int_0^t\big\langle\mu^{I,N}_s,\phi \big\rangle ds+ M^{R,N}_t(\phi) \label{eq:empirical_mu^RN},
	\end{align}
	where for each $e\in{\{S,I,R\}}$, the quantity $M^{e,N}_t$ is a local martingale represented by the following
	
	\begin{align*}
		M^{e,N}_t(\phi) =&\frac{1}{N}\sum_{i=1}^N\int_{0}^{t}\1_{\{E^{i,N}_s=e\}}\dx\phi(\xis)\sigma^e(\xis)dB^i_s \nonumber\\
		&+\frac{1}{N}\sum_{i=1}^N\int_{[0,t]\times\R_+}\big(\1_e(E^{i,N}_{s})-\1_e(E^{i,N}_{s^-})\big)\phi(X^{i,N}_s)\times\\
		&\hspace{3cm}\times\1_{\big\{u\leq K_{\mu^{I,N}_s}(X^{i,N}_s) \1_0(E^{i,N}_{s^-}) + \gamma \1_1(E^{i,N}_{s^-}) \big\}}\bar{Q}^i(ds,du).
	\end{align*}
	
	\begin{remark}
		To avoid confusions, it is worth to notice that we implicitly used three different test function $\phi^S, \phi^I, \phi^R$ for each measure $\mu^{S,N}, \mu^{I,N}, \mu^{R,N}$ in the above system when we perform Itô's calculus.
	\end{remark}
	We know that these local martingles converge to $0$ as $N\to \infty$, and the law of large numbers result established in \cite{Hauray-Pardoux-Vuong22} ensures the convergence of the triple $\big(\mu^{S,N},\mu^{I,N},\mu^{R,N}\big)\in \left(\D\big([0,T],\M(\R^d)\big)\right)^3$ towards $\big(\mu^S,\mu^I,\mu^R\big)\in \left(C\big([0,T],\M(\R^d)\big)\right)^3$, which is the unique solution of the limit system of \eqref{eq:empirical_mu^SN}-\eqref{eq:empirical_mu^RN}.
	
	\smallskip
	Now, if we consider for each epidemiological state the fluctuation process around its mean field limit, namely
	$$\big(\eta^{S,N}, \eta^{I,N}, \eta^{R,N}\big)=\big(\sqrt{N}(\mu^{S,N}-\mu^S),\sqrt{N}(\mu^{I,N}-\mu^I),\sqrt{N}(\mu^{R,N}-\mu^R)\big),$$
	then equation \eqref{evolution_fluctuation_process} becomes the following system:

	\begin{align}
		\big\langle \eta^{S,N}_t,\phi \big\rangle
		=&\big\langle\eta^{S,N}_0,\phi \big\rangle+\int_{0}^{t}\big\langle\eta^{S,N}_s,L^{S,N}_s(\phi)\big\rangle ds+\int_{0}^{t}\big\langle\eta^{N}_s, \langle\mu^S_s, \dx\phi\cdot V^S\rangle\big\rangle ds \nonumber\\
		&-\int_0^t\big\langle \eta^{I,N}_s, \langle\mu^S_s,\phi K\rangle \big\rangle ds + \tilde{M}^{S,N}_t(\phi), \label{eq:fluct_eta^SN}\\
		\big\langle \eta^{I,N}_t,\phi \big\rangle
		=&\big\langle\eta^{I,N}_0,\phi \big\rangle+\int_{0}^{t}\big\langle\eta^{I,N}_s,L^{I,N}_s(\phi)\big\rangle ds+\int_{0}^{t}\big\langle\eta^N_s, \langle\mu^I_s, \dx\phi\cdot V^I\rangle\big\rangle ds \nonumber\\
		&+\int_0^t\big\langle\eta^{S,N}_s, \phi K_{\mu^{I,N}_s}\big\rangle ds + \tilde{M}^{I,N}_t(\phi), \label{eq:fluct_eta^IN}\\
		\big\langle \eta^{R,N}_t,\phi \big\rangle
		=&\big\langle\eta^{R,N}_0,\phi \big\rangle+\int_{0}^{t}\big\langle\eta^{R,N}_s,L^{R,N}_s(\phi)\big\rangle ds+\int_{0}^{t}\big\langle\eta^N_s, \langle\mu^R_s, \dx\phi\cdot V^R\rangle\big\rangle ds \nonumber\\
		&+\gamma\int_0^t\big\langle\eta^{I,N}_s,\phi \big\rangle ds+ \tilde{M}^{R,N}_t(\phi) \label{eq:fluct_eta^RN},
	\end{align}
	where the differential operators $L^{S,N}$, $L^{I,N}$, $L^{R,N}$ are defined by
	\begin{align}
		L^{S,N}_s(\phi)
		=&\frac{1}{2}\tr\big[(\sigma^S{\sigma^S}^{\dag})\dxx\phi\big]+ \dx\phi\cdot V^S_{\mu^N_s}- \phi K_{\mu^{I,N}_s},\\
		L^{I,N}_s(\phi)
		=&\frac{1}{2}\tr\big[(\sigma^I{\sigma^I}^{\dag})\dxx\phi\big]+ \dx\phi\cdot V^I_{\mu^N_s}+ \langle\mu^S_s,\phi K\rangle-\gamma\phi,\\
		L^{R,N}_s(\phi)
		=&\frac{1}{2}\tr\big[(\sigma^R{\sigma^R}^{\dag})\dxx\phi\big]+ \dx\phi\cdot V^R_{\mu^N_s},
	\end{align}
	and the martingale terms $\tilde{M}^{e,N}_t=\sqrt{N}M^{e,N}_t$ for $e\in\{S,I,R\}$.

	\begin{remark}\label{rmk:def_L^SN,L^IN,L^RN}
		
		The first term in the definition of differential operators $L^{S,N}$, $L^{I,N}$, $L^{R,N}$ emerge naturally after renormalizing the difference between the original system \eqref{eq:empirical_mu^SN}-\eqref{eq:empirical_mu^RN} and its limit (there is no linearization here), whereas the other terms represent a part of the linearized terms and the epidemic dynamic.
		
		\smallskip
		We also notice that
		\begin{equation*}
			\int_{E}\eta^N(de)=\int_{E} \big(\eta^{S,N}+\eta^{I,N}+\eta^{R,N}\big)(de).
		\end{equation*}
	\end{remark}
	
	\begin{remark}\label{rmk:martingales_as_distributions}
		We consider the above system as a semimartingale representation of $\eta^{S,N}$, $\eta^{I,N}$, $\eta^{R,N}$ and regard $\tilde{M}^{S,N}$, $\tilde{M}^{I,N}$, $\tilde{M}^{R,N}$ as distributions acting on test functions. More specifically, in the next sections, we will show that they are the distributions in $H^{-(2+2D),D}$. Nevertheless, instead of using the usual notion for the dual product of $\tilde{M}^{e,N}_t$ and function $\phi$, we always write $\tilde{M}^{e,N}_t(\phi)$ to avoid the abuse of notion $\left< \cdot,\cdot \right>$, e.g. when compute the quadratic variations as in \eqref{eq:quadratic_variation_tildeM^SN} below.
	\end{remark}
	
	Before going on, let us give a heuristic description how the limit of the martingale terms should look like. For $e\in \{S,I,R\}$ and any $\phi\in C^2_b(\R^d)$, $\tilde{M}^{e,N}_t(\phi)$ is a real valued martingale with the quadratic variation given by

	
	\begin{align}
		\big\langle \tilde{M}^{S,N}(\phi) \big\rangle _t 
		=&\frac{1}{N}\sum_{i=1}^N\int_{0}^{t}\1_0(E^{i,N}_s)\big(\dx\phi(\xis)\sigma^S(\xis)\big)^2 ds \nonumber\\
		&+\frac{1}{N}\sum_{i=1}^N\int_0^t\1_0(E^{i,N}_{s})\phi(X^{i,N}_s)^2 K_{\mu^{I,N}_s}(X^{i,N}_s) ds, \nonumber\\
		=&\int_0^t \big\langle \mu^{S,N}_s, \big(\dx\phi\sigma^S\big)^2\big\rangle ds +\int_0^t \big\langle \mu^{S,N}_s, \phi^2 K_{\mu^{I,N}_s}\big\rangle ds, \label{eq:quadratic_variation_tildeM^SN}
	\end{align}

	\begin{align}
		\big\langle \tilde{M}^{I,N}(\phi) \big\rangle _t
		=&\frac{1}{N}\sum_{i=1}^N\int_{0}^{t}\1_1(E^{i,N}_s)\big(\dx\phi(\xis)\sigma^I(\xis)\big)^2 ds \nonumber\\
		&+\frac{1}{N}\sum_{i=1}^N\int_0^t\1_0(E^{i,N}_{s})\phi(X^{i,N}_s)^2 K_{\mu^{I,N}_s}(X^{i,N}_s)\nonumber\\
		&+ \frac{1}{N}\sum_{i=1}^N\int_0^t \gamma\1_1(E^{i,N}_{s})\phi(X^{i,N}_s)^2 ds, \nonumber\\
		=&\int_0^t \big\langle \mu^{I,N}_s, \big(\dx\phi\sigma^I\big)^2\big\rangle ds +\int_0^t \big\langle \mu^{S,N}_s, \phi^2 K_{\mu^{I,N}_s}\big\rangle ds + \int_0^t \big\langle \mu^{I,N}_s, \gamma\phi ^2\big\rangle ds, \label{eq:quadratic_variation_tildeM^IN}
	\end{align}
	and
	\begin{align}
		\big\langle \tilde{M}^{R,N}(\phi) \big\rangle _t
		=&\frac{1}{N}\sum_{i=1}^N\int_{0}^{t}\1_2(E^{i,N}_s)\big(\dx\phi(\xis)\sigma^R(\xis)\big)^2 ds \nonumber\\
		&+ \frac{1}{N}\sum_{i=1}^N\int_0^t \gamma\1_1(E^{i,N}_{s})\phi(X^{i,N}_s)^2 ds, \nonumber\\ 
		=&\int_0^t \big\langle \mu^{R,N}_s, \big(\dx\phi\sigma^R\big)^2\big\rangle ds +\int_0^t \big\langle \mu^{I,N}_s, \gamma\phi ^2\big\rangle ds. \label{eq:quadratic_variation_tildeM^RN}
	\end{align}
	
	By the law of large numbers, we can deduce the convergence of the above quadratic variation processes. When $N$ tends to infinity, these processes are determined by the limit measures $\mu^S,\mu^I,\mu^R$ replacing $\mu^{S,N},\mu^{I,N},\mu^{R,N}$ in equations \eqref{eq:quadratic_variation_tildeM^SN}-\eqref{eq:quadratic_variation_tildeM^RN}. And hence if the limit processes $\tilde{M}^{S}$, $\tilde{M}^{I}$, $\tilde{M}^{R}$ (respectively of $\tilde{M}^{S,N}$, $\tilde{M}^{I,N}$, $\tilde{M}^{R,N}$) are continuous martingales with the deterministic quadratic variations, they can be characterized by Gaussian processes.

	\subsection{Main estimates in dual spaces}\label{sec:main_estimates_dual_spaces}
	
	We first establish some estimates for the fluctuations $\eta^{S,N}$, $\eta^{I,N}$, $\eta^{R,N}$ and the martingales $\tilde{M}^{S,N}$, $\tilde{M}^{I,N}$, $\tilde{M}^{R,N}$ with norms in the dual Sobolev spaces $H^{-(1+D), 2D}$ and $H^{-(2+2D), D}$. In our framework, even though the jumps are bounded, the position variables take value in $\R^d$ so the use of weighted Sobolev spaces is necessary. The weights and regularity index of that Sobolev spaces will be identified in the proof and related to the order of moment estimates acquired on the position of individuals.

	\begin{proposition}\label{prop:main_estimates_martingale_parts}
		Under Assumptions \textbf{$\A_1,\, \A_2$}, for any $T>0$ and for each $e\in \{0,1,2\}$, the process $\tilde{M}^{e,N}_t$ is a $H^{-(1+D),2D}$-valued martingale and satisfies
		\begin{equation}
			\sup_{N\ge 1} \E\left[\sup_{t\le T} \big\|\tilde{M}^{e,N}_t\big\|^{2}_{-(1+D),2D} \right] <+\infty.
		\end{equation}
	\end{proposition}
	\begin{proof}
		We give proof for the case of $\tilde{M}^{S,N}_t$. The estimates for $\tilde{M}^{I,N}_t, \tilde{M}^{R,N}_t$ can be obtained by similar arguments.
		
		\smallskip
		Let $(\phi_k)_{k\ge 1}$ be a complete orthonormal basis of $H^{1+D, 2D}$. It suffices to show that
		\begin{equation}
			\sup_{N\ge 1} \sum_{k\ge 1}\E\left[\sup_{t\le T} \big(\tilde{M}^{S,N}_t(\phi_k)\big)^2 \right] <+\infty.
		\end{equation}
		
		Using Doob's inequality and the boundedness of $\sigma, K$, we deduce that
		\begin{align}\label{eq:CLT:proof_uniform_estimtes_MeN_ineq1}
			\sum_{k\ge 1}\E\left[\sup_{t\le T} \big(\tilde{M}^{S,N}_t(\phi_k)\big)^2 \right] \le& C\sum_{k\ge 1}\E\left[\big(\tilde{M}^{S,N}_T(\phi_k)\big)^2 \right]\nonumber\\
			\le& C\sum_{k\ge 1}\E\left[\int_0^T \left\langle \mu^{S,N}, \big(\dx\phi_k\sigma^S\big)^2\right\rangle ds\right]\nonumber\\
			&\hspace{1cm}+C\sum_{k\ge 1}\E\left[\int_0^T \left\langle \mu^{S,N}, \phi_k^2 K_{\mu^{I,N}_s}\right\rangle ds\right]\nonumber\\
			\le& C\sum_{k\ge 1}\E\left[\int_0^T \left\langle \mu^{S,N}, \big(\divv\phi_k\big)^2\right\rangle ds\right]\nonumber\\
			&\hspace{1cm}+C\sum_{k\ge 1}\E\left[\int_0^T \left\langle \mu^{S,N}, \phi_k^2\right\rangle ds\right].
		\end{align}
		
		On the other hand, using the fact that $(X^{i,N},E^{i,N}),\, i=1,\dots,N$ are identically distributed, we have
		\begin{align}\label{eq:CLT:proof_uniform_estimtes_MeN_ineq2}
			r.h.s. =& C\sum_{k\ge 1}\int_0^T\E\left[\frac{1}{N}\sum_{i=1}^N \1_{\{E^{i,N}_s=0\}}\big(\divv\phi_k(X^{i,N}_s)\big)^2 \right]ds\nonumber\\
			&\hspace{1cm}+C\sum_{k\ge 1}\int_0^T\E\left[ \frac{1}{N}\sum_{i=1}^N \1_{\{E^{i,N}_s=0\}}\phi_k^2(X^{i,N}_s)\right]ds\nonumber\\
			\le& C\sum_{k\ge 1}\int_0^T\E\left[\big(\divv\phi_k(X^{1,N}_s)\big)^2 \right]ds+C\sum_{k\ge 1}\int_0^T\E\left[\phi_k^2(X^{1,N}_s)\right]ds.
		\end{align}
		
		Combing \eqref{eq:CLT:proof_uniform_estimtes_MeN_ineq1} and \eqref{eq:CLT:proof_uniform_estimtes_MeN_ineq1} we obtain
		\begin{equation}
			\begin{aligned}
				\sum_{k\ge 1}\E\left[\sup_{t\le T} \big(\tilde{M}^{S,N}_t(\phi_k)\big)^2 \right] \le& C\sum_{k\ge 1}\int_0^T\E\left[\big(\divv\phi_k(X^{1,N}_s)\big)^2 \right]ds\\
				&+C\sum_{k\ge 1}\int_0^T\E\left[\phi_k^2(X^{1,N}_s)\right]ds.
			\end{aligned}
		\end{equation}
		
		Now applying the definition of the linear mappings $\Psi _x, \delta_{x}$ in Lemma \ref{lem:control_linear_forms}, the above inequality can be rewritten as follows
		\begin{align*}
			\sum_{k\ge 1}\E\left[\sup_{t\le T} \big(\tilde{M}^{S,N}_t(\phi_k)\big)^2 \right]
			\le& C\E\left[\int_0^T \big\|\Psi_{X^{1,N}_s}\big\|_{-(1+D),2D}^2 ds\right]\\
			&+C\E\left[\int_0^T \big\|\delta_{X^{1,N}_s}\big\|_{-(1+D),2D}^2 ds\right].
		\end{align*}
		
		Finally, we combine Lemma \ref{lem:control_linear_forms} and Lemma \ref{lem:finite_moment_estimates} to conclude that
		\begin{align*}
			\sup_{N\ge 1} \E\left[\sup_{t\le T} \big\|\tilde{M}^{S,N}_t\big\|^{2}_{-(1+D),2D} \right]\le&\sup_{N\ge 1} \sum_{k\ge 1}\E\left[\sup_{t\le T} \big(\tilde{M}^{S,N}_t(\phi_k)\big)^2 \right] \\
			\le& C\sup_{N\ge 1}\E \left[\sup_{t\le T}\big(1+|X^{1,N}_t|^{4D}\big)\right] < +\infty.
		\end{align*}
	\end{proof}
	
	\begin{proposition}\label{prop:main_estimates_operators_in_dual_space}
		Under Assumptions \textbf{$\A_1,\, \A_2$}, for each $e\in \{0,1,2\}$, for $0<t\le T$ and for every $N$, the operator $L^{e,N}_t$ is a linear continuous mapping from $H^{2 + 2 D, D}$ into $H^{1 + D, 2 D}$ and we have for all $\phi \in H^{2+2 D, D}$,
		\begin{equation}
			\big\|L^{e,N}_t(\phi)\big\|_{1+D, 2 D} \le C_T\big\|\phi\big\|_{2+2 D, D},
		\end{equation}
		where the constant $C_T$ does not depend on $N$ and the randomness.
	\end{proposition}
	
	\begin{proof}
		We recall that
		\begin{align*}
			L^{S,N}_s(\phi)
			=&\frac{1}{2}\tr\big[(\sigma^S{\sigma^S}^{\dag})\dxx\phi\big]+ \dx\phi\cdot V^S_{\mu^N_s}- \phi K_{\mu^{I,N}_s},\\
			L^{I,N}_s(\phi)
			=&\frac{1}{2}\tr\big[(\sigma^I{\sigma^I}^{\dag})\dxx\phi\big]+ \dx\phi\cdot V^I_{\mu^N_s}+ \langle\mu^S_s,\phi K\rangle-\gamma\phi,\\
			L^{R,N}_s(\phi)
			=&\frac{1}{2}\tr\big[(\sigma^R{\sigma^R}^{\dag})\dxx\phi\big]+ \dx\phi\cdot V^R_{\mu^N_s}.
		\end{align*}
		
		Since $\sigma^S, V, K\in C^{1+D}_b$, we easily deduce that
		\begin{equation}\label{eq:sec2:proof:regularity_L:ineq1}
			\big\| L^{S,N}_s(\phi) \big\|_{1+D, 2D}\le \big\|\phi\big\|_{3+D, 2D}\le C\big\|\phi\big\|_{2+2D, D},
		\end{equation}
		where the inequality on the r.h.s follows by the embedding \eqref{embedding:W_to_W} and since the fact that $D\ge 1$.
		
		\smallskip
		The same argument holds true for $L^{R,N}_s(\phi)$.
		
		\smallskip
		In the representation of $L^{I,N}_s(\phi)$, there is an extra term $\left\langle \mu^S_s, \phi K\right\rangle$, which reduces the regularity of the test functions. To treat this tricky term, we start by using the fact that all the derivatives of $K$ up to order $1+D$ are bounded, we can differentiate under the integral sign w.r.t. variable $y$ and obtain the following
		\begin{align*}
			\big\|\left\langle \mu^S_s, \phi K\right\rangle\big\|_{1+D, 2 D}^{2}
			=& \sum_{|k|=0}^{1+D} \int_{\mathbb{R}^{d}} \dfrac{\Big|D^{k}_y\big\langle \mu^S_s, \phi K(\cdot, y)\big\rangle\big|^{2}}{1+|y|^{4 D}} dy \\
			\le& C\int_{\mathbb{R}^{d}} \frac{\big|\left\langle \mu^S_s, \phi\right\rangle\big|^{2}}{1+|y|^{4 D}} dy \\
			\le& C\int_{\mathbb{R}^{d}}|\phi(x)|^2 \mu^S_s(d x) \int_{\R^{d}} \frac{1}{1+|y|^{4 D}} dy.
		\end{align*}
		Using Lemma \ref{lem:control_linear_forms}, we have
		\begin{align*}
			|\phi(x)|^2=&|\delta_x(\phi)|^2\\
			\le& \|\delta_x\|^2_{-(2+2D),D}\|\phi\|^2_{(2+2D),D}\\
			\le& C\big(1+|x|^{2D}\big)\|\phi\|^2_{(2+2D),D}
		\end{align*}
		Hence we deduce that
		\begin{align*}
			\big\|\left\langle \mu^S_s, \phi K\right\rangle\big\|_{1+D, 2 D}^{2}
			\le& C\big\|\phi\big\|_{2+ 2D, D}^{2}  \int_{\R^{d}}\left(1+|x|^{2 D}\right) \mu^S_s(d x) \int_{\R^{d}} \frac{1}{1+|y|^{4 D}} dy\\
			\le& C\big\|\phi\big\|_{2+ 2D, D}^{2} ,
		\end{align*}
		where we get the last inequality by the fact that $4 D>d$ (thus $\int_{\R^{d}} dy/(1+|y|^{4 D})<+\infty$) and $\mu^S_s$ has finite moments of order $2 D$ (by Lemma \ref{lem:finite_moment_estimates}). Again, we can see the essential of weights in the Sobolev spaces in the above proof.
	\end{proof}
	
	\begin{remark}
		In the system \eqref{eq:fluct_eta^SN}-\eqref{eq:fluct_eta^RN}, it remains the terms $$\int_{0}^{t}\big\langle\eta^{N}_s, \langle\mu^e_s, \dx\phi\cdot V^e\rangle\big\rangle ds,\quad e\in\{S,I,R\},$$ which are not involved in the integrals $\int_{0}^{t}\big\langle\eta^{e,N}_s,L^{e,N}_s\big\rangle ds$.
		
		In fact, these terms are created when we linearize the transport terms in \eqref{eq:empirical_mu^SN}-\eqref{eq:empirical_mu^RN}. The functions $\langle\mu^e_s, \dx\phi\cdot V^e\rangle, \; e\in\{S,I,R\}$, they all tested against the distribution $\eta^N_s$. Following the lines in the proof of Proposition \ref{prop:main_estimates_operators_in_dual_space}, we can obtain similar estimates for these functions, i.e. 
		\begin{equation*}
			\big\|\left\langle \mu^e_s, \dx\phi\cdot V^e\right\rangle\big\|_{1+D, 2 D}
			\le C\big\|\phi\big\|_{2+ 2D, D} , \quad e\in\{S,I,R\}.
		\end{equation*}
	\end{remark}
	
	\medskip
	Now we state a proposition concerning the uniform estimate of the fluctuation processes $(\eta^{e,N}_t)_{t\le T}$, $e\in \{0,1,2\}$. This proposition is equivalent to Proposition~\ref{prop:uniform_estimate_etaN}.
	\begin{proposition}\label{prop:main_estimates_etaN_in_dual_space}
		Under Assumptions \textbf{$\A_1,\, \A_2$}, for any $T>0$ and for each $e\in \{0,1,2\}$, the fluctuation process $\eta^{e,N}_t$ belongs to $H^{-(1+D), 2D}$ uniformly in $t$ and $N$, i.e.
		\begin{equation}\label{eq:uniform_estimate_eta}
			\sup_{N\ge 1} \E\left[\sup _{t \le T}\big\|\eta^{e,N}_t\big\|_{-(1+D),2D}^{2}\right]<+\infty.
		\end{equation}
	\end{proposition}
	
	The proof of Proposition \ref{prop:main_estimates_etaN_in_dual_space} is postponed to Section \ref{sec:proof_main_Prop-(1)}.

	\begin{remark}\label{rmk:dual_norm_and_initial_bound} We have in the following some important remarks:
		\begin{itemize}
			\item We have $\big\|\cdot\big\|_{-(2+2 D), D} \le C\big\|\cdot\big\|_{-(1+D), 2 D}$ by the dual embedding \eqref{embedding:Wdual_to_Wdual}. Now combining with Proposition \ref{prop:main_estimates_martingale_parts} and Proposition \ref{prop:main_estimates_etaN_in_dual_space}, we can also ensures that for $e\in \{S,I,R\}$, $\eta^{e,N}_t$ and $\tilde{M}^{e,N}_t$ belong to $H^{-(2+2D), D}$, i.e.
			\begin{align*}
				\sup_{N\ge 1}  \E\left[\sup _{t\le T}\big\|\eta^{e,N}_t\big\|_{-(2+2 D), D}^{2}\right]<& +\infty,\\
				\sup_{N\ge 1}  \E\left[\sup _{t\le T}\big\|\tilde{M}^{e,N}_t\big\|_{-(2+2 D), D}^{2}\right]<& +\infty.
			\end{align*}
			
			In particular, at the initial time, we have $\sup_{N\ge 1} \E\left[\big\|\eta^{e,N}_0\big\|_{-(2+2D), D}^{2}\right]<+\infty$ under the assumptions \textbf{$\A_1$}, \textbf{$\A_2$}.
			
			\smallskip
			\item As a consequence of Proposition \ref{prop:main_estimates_operators_in_dual_space}, we also have the following statement for the adjoint operators: For $e\in \{S,I,R\}$, for every $u\in H^{-(1+D), 2D}$,
			\begin{align}
				\big\|{L^{e,N}_t}^* u\big\|_{-(2+2D), D}^{2} \le& C_T\big\|u\big\|_{-(1+D), 2D}^{2}.
			\end{align}
		\end{itemize}
	\end{remark}
	
	With the above remarks, we can consider the decomposition \eqref{eq:fluct_eta^SN}-\eqref{eq:fluct_eta^RN} as the following adjoint system in $H^{-(2+2D), D}$
	
	\begin{align}
		\eta^{S,N}_t
		=&\eta^{S,N}_0+\int_{0}^{t}{L^{S,N}_s}^*\eta^{S,N}_s ds-\int_{0}^{t}\divv\big(\mu^S_s V^S_{\eta^N_s}\big) ds -\int_0^t \mu^S_s K_{\eta^{I,N}_s} ds + \tilde{M}^{S,N}_t, \label{eq:adjoint:fluct_eta^SN}\\
		\eta^{I,N}_t
		=&\eta^{I,N}_0+\int_{0}^{t}{L^{I,N}_s}^*\eta^{I,N}_s ds-\int_{0}^{t}\divv\big(\mu^I_s V^I_{\eta^N_s}\big) ds +\int_0^t \mu^{S,N}_s K_{\mu^{I,N}_s} ds + \tilde{M}^{I,N}_t, \label{eq:adjoint:fluct_eta^IN}\\
		\eta^{R,N}_t
		=&\eta^{R,N}_0+\int_{0}^{t}{L^{R,N}_s}^*\eta^{R,N}_s ds-\int_{0}^{t}\divv\big(\mu^R_s V^R_{\eta^N_s}\big) ds +\gamma\int_0^t \mu^S_s K_{\eta^{I,N}_s} ds + \tilde{M}^{R,N}_t. \label{eq:adjoint:fluct_eta^RN}
	\end{align}

	%
	
	\subsection{Tightness results}
	
	In the following, we discuss about the benefit of the Hilbert structure of the Sobolev spaces used in this present paper when proving the tightness results. Let us state here the Aldous tightness criterion for Hilbert space valued stochastic processes.
	
	\paragraph{\textbf{Aldous's criterion}}
	(See e.g. \cite{Aldous78}, \cite{Meleard96}) Let $H$ be a separable Hilbert space. A sequence of processes $(X^N)_{N\geq 1}$ in $\mathcal{D}(\mathbb{R}_{+},H)$ defined on the respective filtered probability spaces $(\Omega^N,\mathcal{F}^N,(\mathcal{F}^N_{t})_{t\geq 0},\mathbb{P}^N)$ is tight if it satisfies both the two following conditions:
	
	\smallskip
	\label{condition:A1}
	$\left(\text{A}_{1}\right)$: For every $t\geq 0$ and $\varepsilon>0$, there exists a compact set $K\subset H$ such that $$\sup_{N\geq 1} \mathbb{P}^N\left( X^N_{t}\notin K \right)\leq \varepsilon,$$
	
	\smallskip
	
	\label{condition:A2}
	$\left(\text{A}_{2}\right)$: For every $\varepsilon,\varepsilon_{2}>0$ and $\theta\geq 0$, there exists $\delta_{0}>0$ and an integer $N_{0}$ such that for all $(\mathcal{F}^N_{t})_{t\geq 0}$-stopping time $\tau_{N}\leq \theta$,
	\begin{equation*}
		\sup_{N\geq N_{0}} \sup_{\delta\leq \delta_{0}} \mathbb{P}^N\left( \big\|X^N_{\tau_{N}+\delta} - X^N_{\tau_{N}}\big\|_{H}\geq \varepsilon \right)\leq \varepsilon_{2}.
	\end{equation*}
	
	\smallskip
	
	To check the Aldous criterion, we will use another version of the first condition where \hyperref[condition:A1]{$(\text{A}_{1})$} is replaced by the condition \hyperref[condition:A1']{$(\text{A}_{1}')$} stated below:
	
	\smallskip
	\label{condition:A1'}
	$\left(\text{A}_{1}'\right)$: There exists a Hilbert space $H_{0}$ such that $H_{0}\hookrightarrow_{c} H$ and, for all $t\geq 0$, 
	$$\sup_{N\geq 1} \mathbb{E}^N[\big\| X^N_{t} \big\|^{2}_{H_{0}}]<+\infty,$$ 
	where the notation $\hookrightarrow_{c}$ means that the embedding is compact and $\mathbb{E}^N$ denotes the expectation associated with the probability $\mathbb{P}^N$.
	
	\smallskip
	Indeed, \hyperref[condition:A1]{$(\text{A}_{1})$} is implied by \hyperref[condition:A1']{$(\text{A}_{1}')$} since the embedding $H_{0}\hookrightarrow_{c} H$ is compact, the closed balls in $H_{0}$ are compact in $H$. Combining with the Markov inequality, condition \hyperref[condition:A1]{$(\text{A}_{1})$} is satisfied.
	
	\medskip
	\begin{theorem}\label{thm:tightness_of_M^N}
		The sequences of the laws of $(\tilde{M}^{S,N})_{N \ge 1}$, $(\tilde{M}^{I,N})_{N \ge 1}$, $(\tilde{M}^{R,N})_{N \ge 1}$ are tight in $\D\left([0, T], H^{-(2+2 D), D}\right)$.
	\end{theorem}		
	
	\begin{proof}
		We will only check the two conditions in Aldous's criterion for $\tilde{M}^{S,N}$, the same can be justified for $\tilde{M}^{I,N}$ and $\tilde{M}^{R,N}$.
		
		\smallskip
		Thanks to Proposition \ref{prop:main_estimates_martingale_parts}, condition \hyperref[condition:A1]{$(\text{A}_{1})$} is satisfied with $H_{0}=H^{-(1+D), 2 D}$ and $H=H^{-(2+2 D), D}$ since the embedding $H^{-(1+D), 2 D} \hookrightarrow H^{-(2+2 D), D}$ is compact (see \eqref{embedding:Wdual_to_Wdual}).
		
		Condition \hyperref[condition:A2]{$(\text{A}_{2})$} is obtained as soon as it holds for the trace of the process $\ll \tilde{M}^{S,N} \gg_{t}$,  where $\ll \tilde{M}^{S,N} \gg_{t}$ is the Doob-Meyer process  associated with the martingale $(\tilde{M}^{S,N}_{t})_{t\geq 0}$ and satisfies the following: For any $t> 0$, $\ll \tilde{M}^{S,N} \gg_{t}$ is a linear continuous mapping from $H^{1+D, 2D}$ to $H^{-(1+D), 2D}$ defined for all $\phi$, $\psi$ in $H^{1+D, 2D}$ by
		\begin{equation*}
			\Big\langle {\ll \tilde{M}^{S,N} \gg}_{t} (\phi),\psi \Big\rangle = \int_0^t \big\langle \mu^{S,N}_s, \big(\dx\phi\sigma^S\big)\big(\dx\psi\sigma^S\big)\big\rangle ds +\int_0^t \big\langle \mu^{S,N}_s, \phi \psi K_{\mu^{I,N}_s}\big\rangle ds,
		\end{equation*}
		(See e.g. Rebolledo’s Theorem in \cite{Joffe-Metivier86}).
		
		\smallskip
		Let $T, \varepsilon,\varepsilon_{2}>0$ and let $\tau_{N} \le T$ be a stopping time. For a complete orthonormal basis $(\phi_k)_{k \geq 1}$ in $H^{2+2 D, D}$, we have
		\begin{align*}
			\sup _{N \ge N_{0}} \sup _{\delta \le \delta_0} \;& \Pb\bigg(\left|\tr{\ll \tilde{M}^{S,N} \gg}_{\tau_{N}+\delta}-\tr{\ll \tilde{M}^{S,N} \gg}_{\tau_{N}}\right|>\varepsilon\bigg)\\
			&\le \frac{1}{\varepsilon} \sup _{N \ge N_{0}} \sup _{\delta \le \delta_0}\E\left[\sum_{k \ge 1} \Big\langle{\ll \tilde{M}^{S,N} \gg}_{\tau_{N}+\delta}(\phi_k),\phi_k\Big\rangle-\Big\langle{\ll \tilde{M}^{S,N} \gg}_{\tau_{N}}(\phi_k)\phi_k\Big\rangle\right]\\
			&\le \frac{C}{\varepsilon} \sup _{N \ge N_{0}} \sup _{\delta \le \delta_0}\E\left[\int_{\tau_{N}}^{\tau_{N}+\delta} \left\langle \mu^{S,N}_s, \big\|\Psi _x\big\|_{-(2+2 D), D}^{2} + \big\|\delta_{x}\big\|_{-(2+2 D), D}^{2} \right\rangle ds\right].
		\end{align*}
		
		\smallskip
		At this step, we again use Lemma \ref{lem:control_linear_forms} and Lemma \ref{lem:finite_moment_estimates} to bound the r.h.s.,
		\begin{align*}
			r.h.s.
			&\le \frac{C}{\varepsilon} \sup _{N \ge N_{0}} \sup _{\delta \le \delta_0}\E\left[\int_{\tau_{N}}^{\tau_{N}+\delta} \left\langle \mu^{S,N}_s, \big\|\Psi _x\big\|_{-(1+D), 2D}^{2} + \big\|\delta_{x}\big\|_{-(1+D), 2D}^{2} \right\rangle ds\right] \\
			&\le \frac{C}{\varepsilon} \sup _{N \ge N_{0}} \sup _{\delta \le \delta_0}\E\left[\int_{\tau_{N}}^{\tau_{N}+\delta} \frac{1}{N} \sum_{i=1}^N\left(1+\left|X_{s}^{i, N}\right|^{4 D}\right) ds\right]\\
			&\le \frac{C\delta_0}{\varepsilon}  \sup _{N \ge N_{0}}\E\left[\sup _{s \le T}\left(1+\left|X_{s}^{1, N}\right|^{4 D}\right)\right]\le \varepsilon_{2},
		\end{align*}
		when $\delta_{0}$ is small enough.
		
		\smallskip
		And thus, both the two conditions for tightness are fulfilled.
	\end{proof}
	
	\begin{theorem}\label{thm:tightness_of_eta^N}
		The sequences of the laws of $(\eta^{S,N})_{N \ge 1}$, $(\eta^{I,N})_{N \ge 1}$, $(\eta^{R,N})_{N \ge 1}$ are tight in $\D\left([0, T], H^{-(2+2 D), D}\right)$.
	\end{theorem}

	\begin{proof}
		Proposition \ref{prop:main_estimates_etaN_in_dual_space} implies that condition \hyperref[condition:A1]{$(\text{A}_{1})$} is satisfied with $H_{0}=H^{-(1+D), 2 D}$ and $H=H^{-(2 + 2 D), D}$. Thanks to Rebolledo's Theorem and the proof of Theorem \ref{thm:tightness_of_M^N} for the martingale terms, condition \hyperref[condition:A2]{$(\text{A}_{2})$} for the sequences $(\eta^{e,N})_{N \ge 1}$, $e\in \{S,I,R\}$ are satisfied as soon as they are satisfied for the drift terms. We will check it for the integrals $$\int_{0}^{t} {L^{e,N}_s}^*\big(\eta^{S,N}_s,\eta^{I,N}_s,\eta^{R,N}_s\big)ds,\quad e\in \{S,I,R\},$$ the remaining terms in the adjoint equations \eqref{eq:adjoint:fluct_eta^SN}-\eqref{eq:adjoint:fluct_eta^RN} can be treated in the similar way.
		
		\smallskip
		We now give a proof for instance to $\eta^{S,N}$. Let $T, \varepsilon>0$ and let $\tau_{N} \le T$ be a stopping time. By using Chebyshev's inequality, one can deduce that
		
		\begin{align*}
			\Pb&\left(\bigg\|\int_{0}^{\tau_{N}+\delta} {L^{S,N}_s}^*\eta^{S,N}_s d s-\int_{0}^{\tau_{N}} {L^{S,N}_s}^*\eta^{S,N}_s d s\bigg\|_{-(2+2 D), D} \ge \varepsilon\right)\\
			&\le \frac{1}{\varepsilon^{2}} \E\left[\bigg\|\int_{\tau_{N}}^{\tau_{N}+\delta} {L^{S,N}_s}^*\eta^{S,N}_s d s \bigg\|_{-(2+2 D), D}^{2}\right]\\
			&\le \frac{\delta}{\varepsilon^{2}}  \E\left[\int_{\tau_{N}}^{\tau_{N}+\delta}\big\|{L^{S,N}_s}^*\eta^{S,N}_s\big\|_{-(2+2 D), D}^{2} d s\right].
		\end{align*}
		Let $(\phi_k)_{k \geq 1}$ be a complete orthonormal system in $H^{2+2 D, D}$, we have
		\begin{equation*}
			\big\|{L^{S,N}_s}^*\eta^{S,N}_s\big\|_{-(2+2 D), D}^{2}=\sum_{k \ge 1}\big\langle \eta^{S,N}_s, L^{S,N}_s(\phi_k)\big\rangle^{2}.
		\end{equation*}
		Thus, using Propposition \ref{prop:main_estimates_operators_in_dual_space} we obtain
		\begin{align*}
			r.h.s.
			&\le \frac{\delta}{\varepsilon^{2}}  \E\left[\int_{\tau_{N}}^{\tau_{N}+\delta} \sum_{k \ge 1}\big\langle \eta^{S,N}_s, L^{S,N}_s(\phi_k)\big\rangle^{2} d s\right]\\
			&\le \frac{C\delta}{\varepsilon^{2}} \E\left[\int_{\tau_{N}}^{\tau_{N}+\delta}\big\|\eta^{S,N}_s\big\|_{-(1+D), 2 D}^{2}\right]\\
			&\le \frac{C\delta^2}{\varepsilon^{2}} \E\left[\sup_{s\le T}\big\|\eta^{S,N}_s\big\|_{-(1+D), 2 D}^{2} \right].
		\end{align*}
		
		Now thanks to Proposition \ref{prop:main_estimates_etaN_in_dual_space}, the last expectation is finite and hence, we can find $\delta_0>0$ such that the condition \hyperref[condition:A2]{$(\text{A}_{2})$} is satisfied. The proof for tightness of the laws of $(\eta^{S,N})_{N \ge 1}$ in $\D\big([0,T],H^{-(2+2D), D}\big)$ is completed.
		
	\end{proof}


	\subsection{Proof of Proposition \ref{prop:main_estimates_etaN_in_dual_space}}\label{sec:proof_main_Prop-(1)}	
	
	In this section, we study a semigroup representation of the evolution equation of the fluctuation processes $\eta^{S,N}$, $\eta^{I,N}$, $\eta^{R,N}$. First, we establish the semigroup formalism for the evolution equation of $\eta^{S,N}$, $\eta^{I,N}$, $\eta^{R,N}$ and provide some useful estimates in weighted Sobolev norms related to the regularity of these semigroups. Second, we state a uniform in time estimate for the stochastic convolution with these semigroups. All results obtained in this section are devoted to prove Proposition \ref{prop:main_estimates_etaN_in_dual_space} in Section \ref{sec:main_estimates_dual_spaces}.
	
	\smallskip
	For each epidemiological state $e\in\{S,I,R\}$, we consider the second order differential operator $A^e$ defined by
	\begin{equation}
		A^e:=\frac{1}{2}\dx\cdot\left(\sigma^e{\sigma^e}^{\dag}\dx\right).
	\end{equation}
	The operators $A^e,\, e\in\{S,I,R\}$ are self-adjoint and we have
	\begin{equation}
		\frac{1}{2}\tr[(\sigma^e{\sigma^e}^{\dag})\dxx]=A^e+B\cdot \dx,
	\end{equation}
	where $B=\left(\frac{1}{2}\sum_{i=1}^d \partial_{x_i} \left(\sigma^e{\sigma^e}^{\dag}\right)_{ij}\right)_{1\le j\le d}$.
	
	\smallskip
	Now we introduce a new drift term
	\begin{equation}
		\tilde{V}^e_{\mu}:=V^e_{\mu}+B, \quad e\in\{S,I,R\}.
	\end{equation}
	
	For each $e\in\{S,I,R\}$, we denote by $\big(\T^e_t\big)_{t\ge 0}$ the semigroup generated by $A^e$ on $L^2(\R^d)$.
	First, we show in the following the adjoint equations under the action of these semigroups $\T^S_t, \T^I_t, \T^R_t$.
	
	\begin{lemma}
		For $t\in [0, T]$, the processes $\eta^{S,N}, \eta^{I,N}, \eta^{R,N}$ satisfy the following system:
		\begin{align}
			\eta^{S,N}_t
			=&{\T^S_{t}}^*\eta^{S,N}_0 - \int_{0}^{t}{\T^S_{t-s}}^*\divv\big(\eta^{S,N}_s \tilde{V}^S_{\mu^N_s}\big) ds- \int_{0}^{t}{\T^S_{t-s}}^*\divv\big(\mu^S_s {V}^S_{\eta^N_s}\big) ds \nonumber\\
			&-\int_0^t{\T^S_{t-s}}^*\big(\eta^{S,N}_sK_{\mu^{I,N}_s}\big) ds-\int_0^t {\T^S_{t-s}}^*\big(\mu^S_sK_{\eta^{I,N}_s}\big) ds + \int_0^t{\T^S_{t-s}}^*d\tilde{M}^{S,N}_s, \label{eq:semigroup_fluct_eta^SN}\\
			\eta^{I,N}_t
			=&{\T^I_{t}}^*\eta^{I,N}_0 - \int_{0}^{t}{\T^I_{t-s}}^*\divv\big(\eta^{I,N}_s \tilde{V}^I_{\mu^N_s}\big) ds- \int_{0}^{t}{\T^I_{t-s}}^*\divv\big(\mu^I_s {V}^I_{\eta^N_s}\big) ds \nonumber\\
			&+\int_0^t{\T^I_{t-s}}^*\big(\eta^{S,N}_sK_{\mu^{I,N}_s}\big) ds+\int_0^t {\T^I_{t-s}}^*\big(\mu^S_sK_{\eta^{I,N}_s}\big) ds -\gamma\int_0^t{\T^I_{t-s}}^*\eta^{I,N}_s ds \nonumber\\
			&+\int_0^t{\T^I_{t-s}}^*d\tilde{M}^{I,N}_s, \label{eq:semigroup_fluct_eta^IN}\\
			\eta^{R,N}_t
			=&{\T^R_{t}}^*\eta^{R,N}_0 - \int_{0}^{t}{\T^R_{t-s}}^*\divv\big(\eta^{R,N}_s \tilde{V}^R_{\mu^N_s}\big) ds- \int_{0}^{t}{\T^R_{t-s}}^*\divv\big(\mu^R_s {V}^R_{\eta^N_s}\big) ds \nonumber\\
			&+\gamma\int_0^t{\T^R_{t-s}}^*\eta^{I,N}_s ds + \int_0^t{\T^R_{t-s}}^*d\tilde{M}^{R,N}_s. \label{eq:semigroup_fluct_eta^RN}
		\end{align}
	\end{lemma}
	
	\begin{proof}
		First, we fix $t\in[0,T]$ and $\phi\in C^2\big(\R^d\big)$. Appling Itô's formula to the test function $\psi(s,x)=(\T^S_{t-s}\phi)(x)$, and notice that for all $x\in\R^d$, the mapping $s\mapsto (\T^S_{t-s}\phi)(x)$ is differentiable and 
		\[
		\frac{d}{ds}\T^S_{t-s}\phi(x)= -A^S(\T^S_{t-s}\phi)(x),
		\]
		we can derive the following equation similar to \eqref{eq:fluct_eta^SN},
		\begin{align*}
			\big\langle \eta^{S,N}_t,\phi \big\rangle
			=&\big\langle\eta^{S,N}_0,\T^S_{t}\phi \big\rangle\\
			&+ \int_{0}^{t}\big\langle\eta^{S,N}_s, \dx(\T^S_{t-s}\phi)\cdot \tilde{V}^S_{\mu^N_s}\big\rangle ds+\int_{0}^{t}\big\langle\eta^N_s, \langle\mu^S_s, \dx(\T^S_{t-s}\phi)\cdot {V}^S\rangle\big\rangle ds \nonumber\\
			&-\int_0^t\big\langle\eta^{S,N}_s, (\T^S_{t-s}\phi) K_{\mu^{I,N}_s}\big\rangle ds-\int_0^t\big\langle \eta^{I,N}_s, \langle\mu^S_s,(\T^S_{t-s}\phi) K\rangle \big\rangle ds\\
			&+\int_0^t d\tilde{M}^{S,N}_s(\T^S_{t-s}\phi).
		\end{align*}
	\end{proof}

	\medskip
	
	Before going on, let us provide some useful estimates to control the terms in the system \eqref{eq:semigroup_fluct_eta^SN}-\eqref{eq:semigroup_fluct_eta^RN}. The first one concerns the regularity estimates of the
	semigroups $(\T^e_{t-s})_{e\in \{S,I,R\}}$ in weighted Sobolev spaces $H^{k,\alpha}$ and will be given in the following proposition (see also a more general result in \cite{Hauray-Vuong22}).
	
	\smallskip
	We consider $A$ the second order differential operator given in the divergence form by
	\[A\phi=- \sum_{i,j=1}^{d} \dxi\big(a_{ij}(x)\dxj \phi\big),\] where the coefficients $a_{ij}$ are symmetric, smooth enough (will be precised) and satisfy the uniform ellipticity condition, i.e.
	\[\sum_{i,j=1}^{d} a_{ij}(x)\xi_i\xi_j\ge \lambda |\xi|^2, \quad \forall x, \xi \in \R^d, \]
	for some positive constant $\lambda$. With the above definition, the operator $A$ is a self-adjoint and positive. Let $(\T_t)_{t\ge 0}$ be the semigroup generated by $A$ on $L^2(\R^d)$.
	
	\begin{proposition}\label{prop:estimates_of_semigroups}
		Let $k\ge 0$ and assume that $a_{ij}\in C^{2k+1}_b(\R^d)$. Let $(\T_t)_{t\ge 0}$ be the semigroup generated by $A$. For any $T\ge 0$, there exists a constant $C_T>0$ depends only on $T, d, k, \|a\|_{H^{2k+1, \alpha}}$ such that for any $t\in [0,T]$, the following holds true
		\begin{enumerate}
			\item
			\begin{equation}\label{CLT:eq:weight:semigroup_bound_first_estimate}
				\|\T_t\phi\|_{H^{k, \alpha}}\le C_T\|\phi\|_{H^{k, \alpha}}.
			\end{equation}
			\item
			\begin{equation}\label{CLT:eq:weight:semigroup_bound_second_estimate}
				\|\nabla_x\T_t\phi\|_{H^{k, \alpha}}\le C_T\left(1+\frac{1}{\sqrt{t}}\right)\|\phi\|_{H^{k, \alpha}}.
			\end{equation}
		\end{enumerate}
	\end{proposition}
	
	\medskip
	Another difficulty we need to handle in the system \eqref{eq:semigroup_fluct_eta^SN}-\eqref{eq:semigroup_fluct_eta^RN} is the stochastic convolutions with the semigroups $\T^S_t, \T^I_t, \T^R_t$ , namely
	\[
	\int_0^t{\T^S_{t-s}}^{*}d\tilde{M}^{S,N}_s,\, \int_0^t{\T^S_{t-s}}^{*}d\tilde{M}^{S,N}_s,\, \int_0^t{\T^S_{t-s}}^{*}d\tilde{M}^{S,N}_s.
	\]
	In the following, we provide a first bound for those terms.
	
	\begin{proposition}
		For $0<t\le T$, there exists a positive constant $C_T$ such that
		\begin{align}
			\E\left[\bigg\|\int_0^t{\T^S_{t-s}}^*d\tilde{M}^{S,N}_s\bigg\|^2_{-(1+D),2D}\right] \leq& C_T, \label{eq:bound_stochastic_convolution_semigroup_S}\\
			\E\left[\bigg\|\int_0^t{\T^I_{t-s}}^*d\tilde{M}^{I,N}_s\bigg\|^2_{-(1+D),2D}\right] \leq& C_T, \label{eq:bound_stochastic_convolution_semigroup_I}\\
			\E\left[\bigg\|\int_0^t{\T^R_{t-s}}^*d\tilde{M}^{R,N}_s\bigg\|^2_{-(1+D),2D}\right] \leq& C_T. \label{eq:bound_stochastic_convolution_semigroup_R}
		\end{align}	
	\end{proposition} 
	
	\begin{proof}
		Let $(\phi_k)_{k \geq 1}$ be a complete orthonormal system in $H^{1+ D, 2D}$, we can also using the expression of ${\T^S_{t-s}}^*d\tilde{M}^{S,N}_s$ in $H^{-(1+ D), 2D}$ via this basis, namely
		\begin{align*}
			\E\Bigg[\bigg\|\int_0^t{\T^S_{t-s}}^*d\tilde{M}^{S,N}_s\bigg\|^2_{-(1+D),2D}\Bigg]= \E\Bigg[\int_0^t \sum_{k \ge 1}\big\langle d\tilde{M}^{S,N}_s, \T^S_{t-s}\phi_k\big\rangle^{2} d s\Bigg]
		\end{align*}
		and then have the same estimates follows the lines in the proof of Proposition \ref{prop:main_estimates_martingale_parts}.
	\end{proof}
	
	However the above bounds are not exactly what we need. Instead, we expect to have a uniform in time estimate for the stochastic convolutions with the semigroups by exploiting the independence of the noise terms. Indeed, we can observe that if these terms do not involve a convolution with the semigroups $(\T^e_{t-s})_{e\in \{S,I,R\}}$ then it
	would be a martingale and we can apply the maximal inequalities for a standard martingale, for instance, the Burkholder-Davis-Gundy inequality and obtain the desired bound. On the other hand, even though the convolution with the
	semigroups $(\T^e_{t-s})_{e\in \{S,I,R\}}$ destroys the martingale property, it is still closely related to maximal
	inequalities by the following proposition (See Theorem 2.1 in \cite{Kotelenez84}).
	
	\medskip
	
	\begin{proposition}\label{lem:maximal_inequality_stochastic_convolution_martingale}
		Let $\big(H, \|\cdot\|_H \big)$ be a separable Hilbert space and $\T_{t}$ be a semigroup acting on $H$. We assume the exponential growth condition on $\T_{t}$, $\left\|\T_{t}\right\|_{L(H)} \leq e^{\alpha t}$ for some positive constant $\alpha$. Then, there exists a constant $C>0$ such that for any $H$-valued locally square integrable \text{càdlàg} martingale $M_{t}$,
		$$
		\E\left[ \sup _{0 \leq t \leq T}\left\|\int_{0}^{t} S_{t-s} d M_{s}\right\|^2_H \right] \le C e^{4 \alpha T} \mathbb{E}\left[\big\|M_T\big\|^2_H\right] .
		$$
	\end{proposition}
	
	In \cite{Hamedani-Zangeneh01}, the authors give a generalization for this maximal inequality with $p$-th moment $(0<p<\infty)$ of stochastic convolution integrals.

	\medskip
	Now we are able to prove Proposition \ref{prop:main_estimates_etaN_in_dual_space}.
	
	\begin{proof}[Proof of Proposition \ref{prop:main_estimates_etaN_in_dual_space}.]\label{proof:proof_main_Prop-(1)}
		
		Using the expression in \eqref{eq:semigroup_fluct_eta^SN}, we have
		\begin{equation}\label{eq:CLT:proof_mainProp_ineq1}
			\begin{aligned}
				\big\|\eta^{S,N}_t\big\|_{-(1+D), 2D}
				\le& \big\|{\T^S_{t}}^*\eta^{S,N}_0\big\|_{-(1+D), 2D}\\
				&+\int_{0}^{t}{\big\|\T^S_{t-s}}^*\divv\big(\eta^{S,N}_s \tilde{V}^S_{\mu^N_s}\big)\big\|_{-(1+D), 2D} ds + \int_{0}^{t}{\big\|\T^S_{t-s}}^*\divv\big(\mu^S_s {V}^S_{\eta^N_s}\big)\big\|_{-(1+D), 2D} ds\\
				&+\int_0^t{\big\|\T^S_{t-s}}^*\big(\eta^{S,N}_sK_{\mu^{I,N}_s}\big)\big\|_{-(1+D), 2D} ds+\int_0^t \big\|{\T^S_{t-s}}^*\big(\mu^S_sK_{\eta^{I,N}_s}\big)\big\|_{-(1+D), 2D} ds\\
				&+ \bigg\|\int_0^t{\T^S_{t-s}}^*d\tilde{M}^{S,N}_s\bigg\|_{-(1+D), 2D}.
			\end{aligned}
		\end{equation}

		Let $(\phi_k)_{k \geq 1}$ be a complete orthonormal system in $H^{1+D, 2D}$. Again, we can use the Parseval's identity to represent the dual norms. First we will treat the two terms in the second line on the r.h.s. of \eqref{eq:CLT:proof_mainProp_ineq1}.
		
		\smallskip
		Let us consider the linear mappings $\Phi_1,\, \Phi_2:\, H^{1+D, 2D}\to \R$ defined by
		\begin{align*}
			\Phi_1(\phi_k)=&\big\langle\eta^{S,N}_s, \dx(\T^S_{t-s}\phi_k)\cdot \tilde{V}^S_{\mu^N_s}\big\rangle,\\
			\Phi_2(\phi_k)=&\big\langle\mu^S_s, \dx(\T^S_{t-s}\phi_k)\cdot \tilde{V}^S_{\eta^N_s}\rangle\big\rangle.
		\end{align*}
		
		Using the second inequality in Proposition \ref{prop:estimates_of_semigroups}, we get
		\begin{align*}
			\big|\Phi_1(\phi_k)\big|=&\big|\big\langle\eta^{S,N}_s, \dx(\T^S_{t-s}\phi_k)\cdot \tilde{V}^S_{\mu^N_s}\big\rangle\big|\\
			\le& C\big\|\eta^{S,N}_s\big\|_{-(1+D), 2D} \big\|\dx(\T^S_{t-s}\phi_k)\cdot \tilde{V}^S_{\mu^N_s}\big\|_{1+D, 2D}\\
			\le& C\big\|\eta^{S,N}_s\big\|_{-(1+D), 2D} \big\|\dx(\T^S_{t-s}\phi_k)\big\|_{1+D, 2D}\\
			\le& C_T\left(1+\frac{1}{\sqrt{t-s}}\right)\big\|\eta^{S,N}_s\big\|_{-(1+D), 2D} \big\|\phi_k\big\|_{1+D, 2D}.
		\end{align*}
		We notice that to obtain the third line, we used the assumption that $V\in C^{1+D}_b(\R^d\times\R^d)$. Now by the similar way, we also have
		\begin{equation*}
			\big|\Phi_2(\phi_k)\big|
			\le C_T\left(1+\frac{1}{\sqrt{t-s}}\right)\big\|\mu^S_s\big\|_{-(1+D), 2D} \big\|\phi_k\big\|_{1+D, 2D},
		\end{equation*}
		and using the continuous embedding from $\P(\R^d)$ into $H^{-(1+D), 2D}$, we obtain
		\begin{equation*}
			\big|\Phi_2(\phi_k)\big|
			\le C_T\left(1+\frac{1}{\sqrt{t-s}}\right) \big\|\phi_k\big\|_{1+D, 2D}.
		\end{equation*}
		Hence we deduce that
		\begin{equation}\label{eq:proof_mainProp_bound_1}
			\begin{aligned}
				\int_{0}^{t}{\big\|\T^S_{t-s}}^*&\divv\big(\eta^{S,N}_s \tilde{V}^S_{\mu^N_s}\big)\big\|_{-(1+D), 2D} ds+ \int_{0}^{t}{\big\|\T^S_{t-s}}^*\divv\big(\mu^S_s {V}^S_{\eta^N_s}\big)\big\|_{-(1+D), 2D} ds\\
				\le& \int_{0}^{t}C_T\left(1+\frac{1}{\sqrt{t-s}}\right)\big\|\eta^{S,N}_s\big\|_{-(1+D), 2D} ds + \int_{0}^{t}C_T\left(1+\frac{1}{\sqrt{t-s}}\right) ds.
			\end{aligned}
		\end{equation}
		
		To treat the two terms created by the jumping part in the third line of \eqref{eq:CLT:proof_mainProp_ineq1}, we use the first statement in Proposition \ref{prop:estimates_of_semigroups}. Indeed by the similar way as before, we also obtain the following bounds 
		\begin{align}\label{eq:proof_mainProp_bound_2}
			\begin{split}
				\int_0^t{\big\|\T^S_{t-s}}^*&\big(\eta^{S,N}_sK_{\mu^{I,N}_s}\big)\big\|_{-(1+D), 2D} ds +\int_0^t \big\|{\T^S_{t-s}}^*\big(\mu^S_sK_{\eta^{I,N}_s}\big)\big\|_{-(1+D), 2D} ds\\
				\le& \int_{0}^{t}C_T\big\|\eta^{S,N}_s\big\|_{-(1+D), 2D} ds + \int_{0}^{t}C_T\big\|\mu^S_s\big\|_{-(1+D), 2D} ds\\
				\le& \int_{0}^{t}C_T\big\|\eta^{S,N}_s\big\|_{-(1+D), 2D} ds + C_T.
			\end{split}
		\end{align}
		
		Concerning the last term (the stochastic convolution), we use Proposition \ref{lem:maximal_inequality_stochastic_convolution_martingale}, Jensen's inequality and Proposition \ref{prop:main_estimates_martingale_parts} to deduce the following
		\begin{align}\label{eq:proof_mainProp_bound_3}
			\begin{split}
				\E\left[\sup_{t \le T}\bigg\|\int_0^t{\T^S_{t-s}}^*d\tilde{M}^{S,N}_s\bigg\|_{-(1+D), 2D}\right]\le& \E\left[\sup_{t \le T}\bigg\|\int_0^t{\T^S_{t-s}}^*d\tilde{M}^{S,N}_s\bigg\|^2_{-(1+D), 2D}\right]^{1/2}\\
				\le& C_T\E\left[\big\|\tilde{M}^{S,N}_T\big\|^2_{-(1+D), 2D}\right]^{1/2}\\
				<& +\infty.
			\end{split}
		\end{align}
		
		Now summing up \eqref{eq:proof_mainProp_bound_1}-\eqref{eq:proof_mainProp_bound_3}, we conclude that
		\begin{align*}
			\E\left[\sup_{s\le t}\big\|\eta^{S,N}_s\big\|_{-(1+D), 2D}\right]
			\le& C_T\E\left[\big\|\eta^{S,N}_0\big\|_{-(1+D), 2D}\right]+ \int_{0}^{t}\frac{C_T}{\sqrt{t-s}}\E\left[\big\|\eta^{S,N}_s\big\|_{-(1+D), 2D}\right] ds \\
			& +\int_{0}^{t}C_T\E\left[\big\|\eta^{S,N}_s\big\|_{-(1+D), 2D}\right] ds + C_T.
		\end{align*}
		
		The similar arguments also give us the uniform in time estimates  for $\eta^{I,N}_t$ and $\eta^{R,N}_t$, namely
		\begin{align*}
			\E\left[\sup_{s\le t}\big\|\eta^{I,N}_s\big\|_{-(1+D), 2D}\right]
			\le& C_T\E\left[\big\|\eta^{I,N}_0\big\|_{-(1+D), 2D}\right]+ \int_{0}^{t}\frac{C_T}{\sqrt{t-s}}\E\left[\big\|\eta^{I,N}_s\big\|_{-(1+D), 2D}\right] ds \\
			& +\int_{0}^{t}C_T\E\left[\big\|\eta^{S,N}_s\big\|_{-(1+D), 2D}+\big\|\eta^{I,N}_t\big\|_{-(1+D), 2D}\right] ds + C_T,\\
			\E\left[\sup_{s\le t}\big\|\eta^{R,N}_s\big\|_{-(1+D), 2D}\right]
			\le& C_T\E\left[\big\|\eta^{R,N}_0\big\|_{-(1+D), 2D}\right]+ \int_{0}^{t}\frac{C_T}{\sqrt{t-s}}\E\left[\big\|\eta^{R,N}_s\big\|_{-(1+D), 2D}\right] ds \\
			& +\int_{0}^{t}C_T\E\left[\big\|\eta^{I,N}_s\big\|_{-(1+D), 2D}\right] ds + C_T.
		\end{align*}
		Now combining all the above estimates and denoting \[\varphi(t)=\E\left[\sup_{s\le t}\Big(\big\|\eta^{S,N}_t\big\|_{-(1+D), 2D}+\big\|\eta^{I,N}_t\big\|_{-(1+D), 2D}+\big\|\eta^{R,N}_t\big\|_{-(1+D), 2D}\Big)\right],\]
		we obtain one estimate in type of Gronwall's lemma.
		\begin{align}\label{eq:proof:Granwall-like_estimate}
			\varphi(t)
			\le& C_T\varphi(0)+ C_T\int_{0}^{t}\left(1+\frac{1}{\sqrt{t-s}}\right)\varphi(s) ds + C_T.
		\end{align}
		However it is not straightforward to directly apply Gronwall's lemma to the above estimate. Indeed, we need to do some modifications. By iterating the estimate \eqref{eq:proof:Granwall-like_estimate} we get
		\begin{align}\label{eq:proof:Granwall-like_estimate:ineq1}
			\begin{split}
				\varphi(t)
				\le& \big(C_T\varphi(0)+C_T\big)+ \big(C_T\varphi(0)+C_T\big)C_T\int_{0}^{t}\left(1+\frac{1}{\sqrt{t-s}}\right) ds\\
				&+C_T^2\int_{0}^{t}\int_0^s\left(1+\frac{1}{\sqrt{t-s}}\right)\left(1+\frac{1}{\sqrt{s-r}}\right)\varphi(r) drds\\
				\le& \big(C_T\varphi(0)+C_T\big)\Big(1+C_T\big(T+2\sqrt{T}\big)\Big)\\
				&+C_T^2\int_0^s\varphi(r)\int_{0}^{t}\left(1+\frac{1}{\sqrt{t-s}}\right)\left(1+\frac{1}{\sqrt{s-r}}\right) dsdr,
			\end{split}
		\end{align}
		where we interchanged the order in the integral in the second line.
		
		Now for $r<s<t$, we have
		\begin{align}\label{eq:proof:Granwall-like_estimate:ineq2}
			\begin{split}
				\int_{r}^{t}\left(1+\frac{1}{\sqrt{t-s}}\right)\left(1+\frac{1}{\sqrt{s-r}}\right) ds=&\int_{r}^{t}\left(1+\frac{1}{\sqrt{t-s}}+\frac{1}{\sqrt{s-r}}+\frac{1}{\sqrt{t-s}\sqrt{s-r}}\right) ds\\
				\le& T+2\sqrt{T}+\int_{r}^{t}\frac{ds}{\sqrt{t-s}\sqrt{s-r}}.
			\end{split}
		\end{align}
		By the change of variables $u=s-r,\; v=t-r$ we have
		\begin{align}\label{eq:proof:Granwall-like_estimate:ineq3}
			\begin{split}
				\int_{r}^{t}\frac{ds}{\sqrt{t-s}\sqrt{s-r}} =&\int_0^v \frac{du}{\sqrt{u}\sqrt{v-u}}\\
				\le& \int_0^{v/2} \frac{du}{\sqrt{u}\sqrt{v-u}} + \int_{v/2}^v \frac{du}{\sqrt{u}\sqrt{v-u}}\\
				\le& \frac{1}{\sqrt{v/2}}\int_0^{v/2} \frac{du}{\sqrt{u}} + \frac{1}{\sqrt{v/2}}\int_{v/2}^v \frac{du}{\sqrt{v-u}}\\
				\le& 4.
			\end{split}
		\end{align}
		
		Finally, we combine inequalities \eqref{eq:proof:Granwall-like_estimate:ineq1}, \eqref{eq:proof:Granwall-like_estimate:ineq2} and \eqref{eq:proof:Granwall-like_estimate:ineq3} to obtain an estimate in type of Gronwall's lemma as usual. Using Remark \ref{rmk:dual_norm_and_initial_bound} for the boundedness at the initial time, we complete the proof of Proposition \ref{prop:main_estimates_etaN_in_dual_space}.
		
	\end{proof}

	\section{Characterization of the limit}\label{sec:characterzation_of_the_limit}
	
	The aim of this section is to prove convergence of the sequence of fluctuation processes $(\eta^N)_{N\geq 1}$, where the limit fluctuation processes $\eta$ is the unique solution of a system of SPDEs driven by four inputs: an initial condition and three noises created by the martingale terms $\tilde{M}^{S,N}_t, \tilde{M}^{I,N}_t, \tilde{M}^{R,N}_t$. In Section~\ref{sec:convergence_of_tildeM^N}, we first identify all the noise terms appearing in the limit system. In Section~\ref{sec:convergence_of_eta^N}, we will show that this system uniquely characterizes the limit law and hence complete the proof of the convergence in law of $(\eta^N)_{N\geq 1}$ to $\eta$.

	\subsection{Convergence of $\bigl(\tilde{M}^{S,N},\tilde{M}^{I,N},\tilde{M}^{R,N}\bigr)_{N \geq 1}$}\label{sec:convergence_of_tildeM^N}
	
	Before stating the convergence result of the martingale terms, let us introduce the definition of Gaussian white noises.
	\begin{definition}
		A random distribution $\mathcal{W}$ defined on a probability space $(\Omega, \mathcal{F}, \mathbb{P})$ is called a standard Gaussian white noise on $\mathbb{R}^d$ if the mapping $\varphi \mapsto\left<\mathcal{W}, \varphi\right>$ is linear and continuous from $L^2\left(\mathbb{R}^d\right)$ into $L^2(\Omega)$, and $\left<\mathcal{W}, \varphi\right>$ is a generalized centered Gaussian process satisfying
		$$\mathbb{E}\left[\left<\mathcal{W}, \varphi\right>\left<\mathcal{W}, \phi\right>\right]=\left<\varphi, \phi\right>_{L^2}, \quad \forall\;\varphi, \phi \in L^2(\mathbb{R}^d).$$
		Here $\left<\cdot,\cdot\right>_{L^2}$ denotes a scalar product on $L^2\left(\mathbb{R}^d\right)$.
		
		\smallskip
		Space-time white noise is a Gaussian white noise on $\mathbb{R}_{+} \times \mathbb{R}^d$.
	\end{definition}
	
	\begin{lemma}\label{lem:continuous_law_of_limit_martingales}
		For all $\phi\in H^{-(2+2D), D}$, The limit process $\big(\mathcal{M}^{S}(\phi), \mathcal{M}^{I}, \mathcal{M}^{R}\big)$ of the sequence $\left(\tilde{M}^{S,N}, \tilde{M}^{I,N}, \tilde{M}^{R,N}\right)_{N \geq 1}$ belong a.s. to $\left(C\left(\mathbb{R}_{+}, H^{-(2+2D), D}\right)\right)^{3}$.
	\end{lemma}
	
	\begin{proof}
		We can adapt the proof in the case of real-value \text{càglàg} processes to the \text{càglàg} process taking values in $H^{-(2+2D), D}$. See e.g. Theorem $3.26$ in \cite{Jacod-Shiryaev87}, or \cite{Billingsley99}.
	\end{proof}

	\begin{proposition}\label{prop:convergence_of_martingales}
		The sequence of the martingales $\bigl(\tilde{M}^{S,N},\tilde{M}^{I,N},\tilde{M}^{R,N}\bigr)_{N \geq 1}$ converges in law in $ \left(\D\left(\mathbb{R}_{+}, H^{-(2+2D), D}\right)\right)^{3}$ towards the continuous centered Gaussian process $\left(\mathcal{M}^{S}, \mathcal{M}^{I}, \mathcal{M}^{R}\right)$ with values in $\left(H^{-(2+2 D), D}\right)^3$ defined by: for all $\varphi, \psi, \phi \in H^{2+2D, D}$,
		\begin{align}
			\left\langle\mathcal{M}_{t}^{S}, \varphi\right\rangle &=\int_{0}^{t} \int_{\mathbb{R}^{d}} \sqrt{\mu^S_r(x)}\dx\varphi (x) \sigma^S (x) \mathcal{W}_{1}(d r, d x) \nonumber\\
			&-\int_{0}^{t} \int_{\mathbb{R}^{d}} \sqrt{\mu^S_r(x) \int_{\mathbb{R}^{d}} \mu^I_r(x) K(x, y)d y}  \varphi(x) \mathcal{W}_{2}(d r, d x), \label{eq:limit_M^S}\\
			\left\langle\mathcal{M}_{t}^{I}, \psi\right\rangle &=\int_{0}^{t} \int_{\mathbb{R}^{d}} \sqrt{\mu^I_r(x)}\dx\psi (x) \sigma^I (x) \mathcal{W}_{1}(d r, d x) \nonumber\\
			&+\int_{0}^{t} \int_{\mathbb{R}^{d}} \sqrt{\mu^S_r(x) \int_{\mathbb{R}^{d}} \mu^I_r(x) K(x, y)d y}  \psi(x) \mathcal{W}_{2}(d r, d x) \nonumber\\
			&-\int_{0}^{t} \int_{\mathbb{R}^{d}} \psi(x) \sqrt{\mu^I_r(x)} \mathcal{W}_{3}(d r, d x), \label{eq:limit_M^I}\\
			\left\langle\mathcal{M}_{t}^{R}, \phi\right\rangle &=\int_{0}^{t} \int_{\mathbb{R}^{d}} \sqrt{\mu^R_r(x)}\dx\phi (x) \sigma^R (x) \mathcal{W}_{1}(d r, d x) \nonumber\\
			&+\int_{0}^{t} \int_{\mathbb{R}^{d}} \phi(x) \sqrt{\gamma \mu^I_r(x)} \mathcal{W}_{3}(d r, d x), \label{eq:limit_M^R}
		\end{align}
		where $\mathcal{W}_{1}, \mathcal{W}_{2}, \mathcal{W}_{3}$ are independent standard space-time white noises.
	\end{proposition}
	
	\begin{proof}
		In the previous section, we proved that the sequence $\big(\tilde{M}^{S,N}, \tilde{M}^{S,N}, \tilde{M}^{S,N}\big)_{N \geq 1}$ is tight in $\left(\D\big(\mathbb{R}_{+}, H^{-(2+2D), D}\big)\right)^{3}$. Hence, according to Prokhorov's Theorem, there exists a subsequence (still denoted by $\big(\tilde{M}^{S,N}, \tilde{M}^{S,N}, \tilde{M}^{S,N}\big)_{N \geq 1}$), which converges in law in $\left(\D\left(\mathbb{R}_{+}, H^{-(2+2D), D}\right)\right)^{3}$ towards $\left(\mathcal{M}^{S}, \mathcal{M}^{I}, \mathcal{M}^{R}\right)$.
		
		\smallskip
		For all $\phi_1,\phi_2,\phi_3 \in H^{2+2D, D}$, by Lemma \ref{lem:continuous_law_of_limit_martingales}, we know that $\mathcal{M}^{S}(\phi_1), \mathcal{M}^{I}(\phi_2), \mathcal{M}^{R}(\phi_3)$ are continuous martingales and thus for any $a_1, a_2, a_3\in \R$, $a_1\mathcal{M}^{S}(\phi_1)+ a_2\mathcal{M}^{I}(\phi_2)+ a_3\mathcal{M}^{R}(\phi_3)$ is also a continuous martingale. Now, we will show that the centered, continuous martingale $\big(\mathcal{M}^{S}(\phi_1), \mathcal{M}^{I}(\phi_2), \mathcal{M}^{R}(\phi_3)\big)$ is a Gaussian process and satisfies \eqref{eq:limit_M^S}-\eqref{eq:limit_M^R}.
		
		\smallskip
		Indeed, let us identify the limit. The LLN result implies that $\big(\mu^{S,N},\mu^{I,N},\mu^{R,N}\big)$ converges in $\left(\D\big([0,T],\M(\R^d)\big)\right)^3$ towards $\big(\mu^S,\mu^I,\mu^R\big)$, which is the unique solution of the limit system of \eqref{eq:empirical_mu^SN}-\eqref{eq:empirical_mu^RN}, and we have
		\begin{align*}
			\big\langle \mathcal{M}^{S}(\phi) \big\rangle _t 
			=&\int_0^t \big\langle \mu^{S}_s, \big(\dx\phi\sigma^S\big)^2\big\rangle ds +\int_0^t \big\langle \mu^{S}_s, \phi^2 K_{\mu^{I}_s}\big\rangle ds, \\
			\big\langle \mathcal{M}^{I}(\phi) \big\rangle _t 
			=&\int_0^t \big\langle \mu^{I}_s, \big(\dx\phi\sigma^I\big)^2\big\rangle ds +\int_0^t \big\langle \mu^{S}_s, \phi^2 K_{\mu^{I}_s}\big\rangle ds + \int_0^t \big\langle \mu^{I}_s, \gamma\phi ^2\big\rangle ds, \\
			\big\langle \mathcal{M}^{R}(\phi) \big\rangle _t 
			=&\int_0^t \big\langle \mu^{R}_s, \big(\dx\phi\sigma^R\big)^2\big\rangle ds +\int_0^t \big\langle \mu^{I}_s, \gamma\phi ^2\big\rangle ds. 
		\end{align*}
		It turns out that $\left\langle a_1\mathcal{M}^{S}(\phi_1)+ a_2\mathcal{M}^{I}(\phi_2)+ a_3\mathcal{M}^{R}(\phi_3)\right\rangle_t$ is a continuous martingale with a deterministic quadratic variation, so it is characterized as a Gaussian process determined by \eqref{eq:limit_M^S}-\eqref{eq:limit_M^R}.
	\end{proof}

	\subsection{Convergence of $\bigl(\eta^{S,N},\eta^{I,N},\eta^{R,N}\bigr)_{N \geq 1}$}\label{sec:convergence_of_eta^N}
	
	We now prove convergence of the sequence $\bigl(\eta^{S,N},\eta^{I,N},\eta^{R,N}\bigr)_{N \geq 1}$ and give a characterization of the limit processes as solution of an equation in $H^{-(4+2 D), D}$. We consider the Hilbert semimartingale decomposition \eqref{eq:adjoint:fluct_eta^SN}-\eqref{eq:adjoint:fluct_eta^RN} of $\bigl(\eta^{S,N},\eta^{I,N},\eta^{R,N}\bigr)$, and we will find a semimartingale decomposition for the limit values, denoted by $\bigl(\eta^S,\eta^I,\eta^R\bigr)$. The difficulty is to close this limit decomposition, i.e. to find a good space in which to immerse the limit process and which allows to give a sense to the limit drift terms. We have seen that the processes $\eta^{S,N},\eta^{I,N},\eta^{R,N}$ belong uniformly to $H^{-(1+D), 2 D}$ and are tight in $H^{-(2+2 D), D}$. We also know that the limit processes $\eta^S,\eta^I,\eta^R$ are in $H^{-(2+2 D), D}$. But to identify the limit in the drift terms, we need to work in a large space that is $H^{-(4+2 D), D}$. And this will be possible if we assume more regularity on the coefficients $\sigma$ and $b$.
	
	\medskip
	We now introduce the following limit operators $L^S_s, L^I_s, L^R_s$ of the linear operators $L^{S,N}$, $L^{I,N}$, $L^{R,N}$, defined by
	
	\begin{align}
		L^{S}_s(\phi)
		=&\frac{1}{2}\tr\big[(\sigma^S{\sigma^S}^{\dag})\dxx\phi\big]+ \dx\phi\cdot V^S_{\mu_s}- \phi K_{\mu^{I}_s},\\
		L^{I}_s(\phi)
		=&\frac{1}{2}\tr\big[(\sigma^I{\sigma^I}^{\dag})\dxx\phi\big]+ \dx\phi\cdot V^I_{\mu_s}+ \langle\mu^S_s,\phi K\rangle-\gamma\phi,\\
		L^{R}_s(\phi)
		=&\frac{1}{2}\tr\big[(\sigma^R{\sigma^R}^{\dag})\dxx\phi\big]+ \dx\phi\cdot V^R_{\mu_s}.
	\end{align}

	Under the Assumption \hyperref[assumption_A3]{\textbf{$\A_3$}}, and follows the lines in the proof of Proposition \ref{prop:main_estimates_operators_in_dual_space}, we can also prove the following lemma.
	
	\begin{lemma}\label{lem:limit_linear_operators}
		For $e\in \{S,I,R\}$, for every $N$ and any $t\le T$, the operators $L^e_s, L^{e,N}_s: H^{4+2 D, D} \rightarrow H^{2+2D, D}$ are linear, continuous and satisfies
		\begin{align}
			\big\|L^{e,N}_t(\phi)\big\|_{2+2D, D} \le& C_T\big\|\phi\big\|_{4+2 D, D},\\
			\big\|L^{e}_t(\phi)\big\|_{2+2D, D} \le& C_T\big\|\phi\big\|_{4+2 D, D}.
		\end{align}
		where the constant $C_T$ does not depend on $N$ and the randomness.
	\end{lemma}
	
	Now by the trivial embedding $H^{-(2+2 D), D} \hookrightarrow$ $H^{-(4+2 D), D}$, the sequence $\bigl(\eta^{S,N},\eta^{I,N},\eta^{R,N}\bigr)_{N \geq 1}$ also converges to $\bigl(\eta^S,\eta^I,\eta^R\bigr)$ in $\left(C\left([0, T], H^{-(4+2 D), D}\right)\right)^3$. This result is stated by the following theorem.
	
	\begin{theorem}
		Under Assumptions $\A_1, \A_2, \A_3$, the sequence $\bigl(\eta^{S,N},\eta^{I,N},\eta^{R,N}\bigr)_{N \geq 1}$ converges in law in $\left(\D\left([0, T], H^{-(2+2 D), D}\right)\right)^3$ to a process $\bigl(\eta^S,\eta^I,\eta^R\bigr)$ which solves the following equation
		
		\begin{align}
			&\eta^{S}_t-\eta^{S}_0-\int_{0}^{t}{L^{S}_s}^*\eta^{S}_s ds+\int_{0}^{t}\divv\big(\mu^S_s V^S_{\eta^S_s+\eta^I_s+\eta^R_s}\big) ds +\int_0^t \mu^S_s K_{\eta^{I}_s} ds=\mathcal{M}^{S}_t, \label{eq:adjoint:limit_fluct_eta^S}\\
			&\eta^{I}_t-\eta^{I}_0-\int_{0}^{t}{L^{I}_s}^*\eta^{I}_s ds+\int_{0}^{t}\divv\big(\mu^I_s V^I_{\eta^S_s+\eta^I_s+\eta^R_s}\big) ds -\int_0^t \mu^{S}_s K_{\mu^{I}_s} ds=\mathcal{M}^{I}_t, \label{eq:adjoint:limit_fluct_eta^I}\\
			&\eta^{R}_t-\eta^{R}_0-\int_{0}^{t}{L^{R}_s}^*\eta^{R}_s ds+\int_{0}^{t}\divv\big(\mu^R_s V^R_{\eta^S_s+\eta^I_s+\eta^R_s}\big) ds -\gamma\int_0^t \mu^S_s K_{\eta^{I}_s} ds=\mathcal{M}^{R}_t, \label{eq:adjoint:limit_fluct_eta^R}
		\end{align}
		where $\mathcal{M}^{S}, \mathcal{M}^{I}, \mathcal{M}^{R}$ are the Gaussian processes defined in Proposition~\ref{prop:convergence_of_martingales}.
	\end{theorem}
	
	\begin{proof}
		Since the sequence of the martingale terms $\left(\tilde{M}^{S,N}, \tilde{M}^{S,N}, \tilde{M}^{S,N}\right)_{N \geq 1}$ converges in law in $\left(\D\left(\mathbb{R}_{+}, H^{-(2+2D), D}\right)\right)^{3}$ to the Gaussian vector process $\left(\mathcal{M}^{S}, \mathcal{M}^{I}, \mathcal{M}^{R}\right)$ defined in Proposition \ref{prop:convergence_of_martingales}, thus to prove that the limit processes satisfies the system \eqref{eq:adjoint:limit_fluct_eta^S}-\eqref{eq:adjoint:limit_fluct_eta^R}, it suffices to show that
		
		\begin{align*}
			&\eta^{S,N}_t-\eta^{S,N}_0-\int_{0}^{t}{L^{S,N}_s}^*\eta^{S,N}_s ds+\int_{0}^{t}\divv\big(\mu^S_s V^S_{\eta^N_s}\big) ds +\int_0^t \mu^S_s K_{\eta^{I,N}_s} ds,\\
			&\eta^{I,N}_t-\eta^{I,N}_0-\int_{0}^{t}{L^{I,N}_s}^*\eta^{I,N}_s ds+\int_{0}^{t}\divv\big(\mu^I_s V^I_{\eta^N_s}\big) ds -\int_0^t \mu^{S,N}_s K_{\mu^{I,N}_s} ds, \\
			&\eta^{R,N}_t-\eta^{R,N}_0-\int_{0}^{t}{L^{R,N}_s}^*\eta^{R,N}_s ds+\int_{0}^{t}\divv\big(\mu^R_s V^R_{\eta^N_s}\big) ds -\gamma\int_0^t \mu^S_s K_{\eta^{I,N}_s} ds
		\end{align*}
		converges in law to
		\begin{align*}
			&\eta^{S}_t-\eta^{S}_0-\int_{0}^{t}{L^{S}_s}^*\eta^{S}_s ds+\int_{0}^{t}\divv\big(\mu^S_s V^S_{\eta^S_s+\eta^I_s+\eta^R_s}\big) ds +\int_0^t \mu^S_s K_{\eta^{I}_s} ds,\\
			&\eta^{I}_t-\eta^{I}_0-\int_{0}^{t}{L^{I}_s}^*\eta^{I}_s ds+\int_{0}^{t}\divv\big(\mu^I_s V^I_{\eta^S_s+\eta^I_s+\eta^R_s}\big) ds -\int_0^t \mu^{S}_s K_{\mu^{I}_s} ds,\\
			&\eta^{R}_t-\eta^{R}_0-\int_{0}^{t}{L^{R}_s}^*\eta^{R}_s ds+\int_{0}^{t}\divv\big(\mu^R_s V^R_{\eta^S_s+\eta^I_s+\eta^R_s}\big) ds -\gamma\int_0^t \mu^S_s K_{\eta^{I}_s} ds,
		\end{align*}
		when $N$ tends to $\infty$. By Lemma \ref{lem:limit_linear_operators}, the integrals $\int_{0}^{t}{L^{S}_s}^*\eta^{S}_s ds,\, \int_{0}^{t}{L^{I}_s}^*\eta^{I}_s ds,\, \int_{0}^{t}{L^{R}_s}^*\eta^{R}_s ds$ and the remaining drift terms make sense in $H^{-(4+2D), D}$. Now, for any $\phi \in H^{-(4+2D), D}$, let us introduce linear vector function
		$F^{\phi}=\left(F^{S,\phi},\, F^{I,\phi},\, F^{R,\phi}\right)$ from $\left(\D\left([0,T], H^{-(2+2D), D}\right)\right)^3$ into $\R^3$ defined by
		\begin{align*}
			F^{S,\phi}_t(u)=&\big\langle u_t,\phi \big\rangle-\big\langle u_0,\phi \big\rangle-\int_{0}^{t}\big\langle u_s,L^{S}_s(\phi)\big\rangle ds\\
			&-\int_{0}^{t}\big\langle (u_s+v_s+w_s), \langle\mu^S_s, \dx\phi\cdot V^S\rangle\big\rangle ds +\int_0^t\big\langle v_s, \langle\mu^S_s,\phi K\rangle \big\rangle ds,\\
			F^{I,\phi}_t(v)=&\big\langle v_t,\phi \big\rangle-\big\langle v_0,\phi \big\rangle-\int_{0}^{t}\big\langle v_s,L^{I}_s(\phi)\big\rangle ds\\
			&-\int_{0}^{t}\big\langle (u_s+v_s+w_s), \langle\mu^I_s, \dx\phi\cdot V^I\rangle\big\rangle ds-\int_0^t\big\langle u_s, \phi K_{\mu^{I}_s}\big\rangle ds,\\
			F^{R,\phi}_t(w)=&\big\langle w_t,\phi \big\rangle-\big\langle w_0,\phi \big\rangle-\int_{0}^{t}\big\langle w_s,L^{R}_s(\phi)\big\rangle ds\\
			&-\int_{0}^{t}\big\langle (u_s+v_s+w_s), \langle\mu^R_s, \dx\phi\cdot V^R\rangle\big\rangle ds-\gamma\int_0^t\big\langle v_s,\phi \big\rangle ds.
		\end{align*}
		
		The function $F^{\phi}$ is continuous and thus, the sequence $\left(F^{\phi}(\eta^{S,N},\eta^{I,N},\eta^{R,N})\right)_{N\ge 1}$ converges in law to $\left(F^{\phi}(\eta^{S},\eta^{I},\eta^{R})\right)$ since the sequence $\bigl(\eta^{S,N},\eta^{I,N},\eta^{R,N}\bigr)_{N \geq 1}$ converges in law to $\bigl(\eta^S,\eta^I,\eta^R\bigr)$ by the tightness result \ref{thm:tightness_of_eta^N}.
		
		\smallskip
		Now it remains to show that $\int_{0}^{t}\big\langle \eta^{S,N}_s,L^{S,N}_s(\phi)-L^{S}_s(\phi)\big\rangle ds$ (and the analogues for $\eta^{I,N}_s,\, \eta^{R,N}_s$) tends to $0$ when $N$ tends to $\infty$. We will prove that it tends to $0$ in $L^1$. Indeed, by Cauchy-Schwartz's inequality, we deduce that
		\begin{align*}
			\E\bigg[\int_0^t \Big|\big\langle & \eta^{S,N}_s,L^{S,N}_s(\phi)-L^{S}_s(\phi)\big\rangle\Big|ds \bigg]\\
			\le& \E\left[\int_0^t\big\| \eta^{S,N}_s\big\|_{-(2+2D,D)^2}\big\|L^{S,N}_s(\phi)-L^{S}_s(\phi)\big\|_{2+2D, D}^2 ds \right]\\
			\le& \int_0^t\E\left[\big\| \eta^{S,N}_s\big\|_{-(2+2D,D)}^2\right]^{1/2}\E\left[\big\|L^{S,N}_s(\phi)-L^{S}_s(\phi)\big\|_{2+2D, D}^2 \right]^{1/2} ds\\
			\le& C\int_0^t\E\left[\big\|L^{S,N}_s(\phi)-L^{S}_s(\phi)\big\|_{2+2D, D}^2 \right]^{1/2} ds,
		\end{align*}
		where we used Proposition \ref{prop:main_estimates_etaN_in_dual_space} and Remark \ref{rmk:dual_norm_and_initial_bound} to obtain the last inequality.
		
		\smallskip
		Following the lines in the proof of Proposition \ref{prop:main_estimates_operators_in_dual_space} and the LLN result $\mu^{e,N}\rightarrow \mu^{e}$, $e\in \{S,I,R\}$, we can also prove that 
		$\big\|L^{S,N}_s(\phi)-L^{S}_s(\phi)\big\|_{2+2D, D}$ tends to $0$ as $N$ tends to $\infty$, and thus complete the proof.
		
		\smallskip
		Noticing that to compute $\big\|L^{S,N}_s(\phi)-L^{S}_s(\phi)\big\|_{2+2D, D}$, we used the additional assumption on $\sigma, V, K$, and once we do the analysis for the term $\frac{1}{2}\tr\big[(\sigma^e{\sigma^e}^{\dag})\dxx\phi\big]$ in $L^{S,N}_s(\phi), L^{S}_s(\phi)$, it will require the regularity order $4+2D$ instead of $2+2D$ as in the inequality \eqref{eq:sec2:proof:regularity_L:ineq1}. Thus, the equations \eqref{eq:adjoint:limit_fluct_eta^S}-\eqref{eq:adjoint:limit_fluct_eta^R} are regarded as the equations in the space $H^{-(4+2D), D}$, while $\eta^S,\eta^I,\eta^R$ are known to take values in the smaller space $H^{-(2+2D), D}$.
	\end{proof}
	
	\medskip
	In order to complete the proof of convergence of the sequence $\bigl(\eta^{S,N},\eta^{I,N},\eta^{R,N}\bigr)_{N \geq 1}$, it remains to prove the uniqueness of solutions to the system \eqref{eq:adjoint:limit_fluct_eta^S}-\eqref{eq:adjoint:limit_fluct_eta^R}.
	
	\begin{proposition}
		For any initial condition $\eta^S_0,\eta^I_0,\eta^R_0$ with values in $H^{-(2+2D), D}$, the system \eqref{eq:adjoint:limit_fluct_eta^S}-\eqref{eq:adjoint:limit_fluct_eta^R} has at most one solution with paths in $\left(\D\left([0,T], H^{-(2+2D), D}\right)\right)^3$.
	\end{proposition}
	
	Since the equations \eqref{eq:adjoint:limit_fluct_eta^S}-\eqref{eq:adjoint:limit_fluct_eta^R} are linear, the standard argument is to take two solutions of this system with the same initial condition and paths in $\left(\D\left([0,T], H^{-(2+2D), D}\right)\right)^3$. Considering an orthonormal basis of $H^{4+2 D, D}$, we can prove the uniqueness of solutions to this system \eqref{eq:adjoint:limit_fluct_eta^S}-\eqref{eq:adjoint:limit_fluct_eta^R} in $\left(\D\left([0,T], H^{-(2+2D), D}\right)\right)^3$. For instance, we can follows the same argument as the proof of uniqueness in \cite{Kurtz-Xiong04}.

\end{document}